\documentclass{amsart}

\usepackage{amssymb,amsmath,epsfig,graphics,latexsym}
\usepackage[english]{babel}
\usepackage[all]{xy}
\newtheorem{theorem}{Theorem}[section]
\newtheorem{proposition}[theorem]{Proposition}
\newtheorem{corollary}[theorem]{Corollary}
\newtheorem{lemma}[theorem]{Lemma}
\newtheorem{claim}[theorem]{Claim}

\theoremstyle{definition}
\newtheorem{definition}[theorem]{Definition}

\theoremstyle{remark}
\newtheorem{remark}[theorem]{Remark}

\newtheorem{question}{Question}
\numberwithin{equation}{section}
\renewcommand{\t}{ \tilde}
\renewcommand{\b}{ \partial}
\newcommand{\Z}{\bf Z}
\newcommand{\R}{\bf R}
\newcommand{\N}{\bf N}

\newcommand{\Hi}{\bf H}

\newcommand{\C}{\bf C}

\renewcommand{\S}{\bf S}
\renewcommand{\l}{\langle}
\renewcommand{\r}{\rangle}

\renewcommand{\int}[1]{\stackrel{\circ}{#1}}
\renewcommand{\o}{\overline}

\newcommand{\co}{\colon\thinspace}
\renewcommand{\epsilon}{\varepsilon}
\newcommand{\cal}{\mathcal}

\usepackage{times}

\begin{document}
\sloppy

\title[]{Non-zero degree maps between closed orientable three-manifolds}

\author{Pierre Derbez}
\address{Laboratoire d'Analyse, Topologie et Probabilit\'es, UMR 6632, 
Centre de Math\'ematiques et d'Informatique, Universit\'e  Aix-Marseille I,
Technopole de Chateau-Gombert, 
39, rue Fr\'ed\'eric Joliot-Curie -
 13453 Marseille Cedex 13}
\email{derbez@cmi.univ-mrs.fr}



\maketitle

\begin{abstract}
This paper adresses the following problem: Given a closed orientable three-manifold $M$, are there at most finitely many closed orientable three-manifolds 1-dominated by $M$? We solve this question for the class of closed orientable graph manifolds. More presisely the main result of this paper asserts that  any closed orientable graph manifold 1-dominates at most finitely many orientable closed three-manifolds satisfying the Poincar\'e-Thurston Geometrization Conjecture. To prove this result we state a more general theorem for Haken manifolds which says that  any closed orientable three-manifold $M$ 1-dominates at most finitely many Haken manifolds whose Gromov simplicial volume is sufficiently close to that of $M$.   
\end{abstract}
\section{Introduction}
\subsection{Statement of the general problem}
We deal here with non-zero degree maps between closed orientable 3-manifolds. Recall that a 3-manifold is termed \emph{geometric} if it admits one of the eight uniform geometries classified by W. P.  Thurston.  Denote by ${\mathcal G}$   the set of closed geometric and Haken manifolds union  the connected sums of such manifolds. Note that the Poincar\'e-Thurston Geometrization Conjecture asserts that ${\mathcal G}$ represents all closed orientable 3-manifolds. Thus a 3-manifold of ${\mathcal G}$ will be termed a \emph{Poincar\'e-Thurston 3-manifold}.  According to \cite{BW}, given two closed orientable 3-manifolds $M$, $N$,   we say that  $M$ $d$-dominates  $N$ ($M\geq_{(d)}N$) if there is a map $f\co M\to N$ of degree $d\not=0$. A motivation for studying nonzero degree maps comes from the observation that they seem to give a way to measure the topological complexity of 3-manifolds and of knots in ${\bf S}^3$. For instance Y. Rong proved in \cite{Ro2} that degree one maps define a partial order on the set ${\mathcal G}$, up to homotopy equivalence. In the same way one can define a partial order on the set ${\mathcal K}$ of knots in ${\bf S}^3$, up to knots equivalence. Given two knots $K$ and $K'$ in  ${\mathcal K}$ we say that $K$ 1-dominates $K'$ if the complement $E_K$ of $K$ properly 1-dominates $E_{K'}$. Then it follows from \cite{Wa} combined with the fact that knots in ${\bf S}^3$ are determined by their complement, see \cite{GL}, that $({\mathcal K},\geq_{(1)})$ is a partially ordered set (a \emph{poset}). This paper adresses the following question which is closely related to the partial order induced by degree one maps (see also Kirby's Problem List \cite[Problem 3.100]{K}):
\begin{question}\label{lf}  Given a closed orientable 3-manifold  $M$,  are there at most finitely many  3-manifolds $N$ in ${\mathcal G}$ (up to homeomorphism) 1-dominated by $M$?
\end{question}
Note that in this question the targets are 3-manifolds of ${\mathcal G}$ because of the Poincar\'e Conjecture. Indeed if there is a fake 3-sphere $K$  then one can get infinitely many reducible homotopy 3-spheres by doing connected sums of finitely many copies of $K$ and since there always exists a degree one map from a closed orientable 3-manifold $M$ to a homotopy 3-sphere we have to exclude this kind of 3-manifolds. On the other hand, in Question 1, we consider always degree one maps to avoid some easy counter examples. For instance for any spherical Lens space $L(p,q)$ there always exists a nonzero degree map (actually a finite covering) from the 3-sphere ${\bf S}^3$ to $L(p,q)$.
\subsection{The main result}
In this paper we solve Question \ref{lf} when the domain $M$ is a closed orientable graph manifold. More precisely our main result states as follows.
\begin{theorem}\label{poset}
Any closed orientable graph manifold 1-dominates at most finitely many closed orientable Poincar\'e-Thurston 3-manifolds. 
\end{theorem}
This result comes from a more general theorem which gives an affirmative answer to Question 1 when the targets are closed Haken manifolds whose Gromov simplicial volume, denoted by ${\rm Vol}(.)$, is sufficiently close to that of the domain $M$. More precisely:

\begin{theorem}\label{dom} For any closed orientable 3-manifold   $M$  there exists   a constant $c\in(0,1)$, which depends only on $M$, such  that $M$  1-dominates at most finitely many  closed Haken manifolds $N$ satisfying ${\rm Vol}(N)\geq(1-c){\rm Vol}(M)$.
\end{theorem}
Recall that there are many important results related to Question \ref{lf} obtained when the targets are restricted. More precisely the known answers can be summurized as follows.
\begin{theorem}[\cite{H-LWZ},\cite{S},\cite{WZ},\cite{Re},\cite{Ro1}]\label{known} Any closed orientable 3-manifold  1-dominates at most finitely many orientable closed geometric 3-manifolds.
\end{theorem}
Notice that in some cases the degree of the maps need not to be bounded. This is true in particular when the targets admit a hyperbolic or an ${\bf H}^2\times{\R}$-structure. Thus a useful consequence of the proof Theorem \ref{known} is the following result.
\begin{corollary}[\cite{S},\cite{WZ}]\label{relative}
Any  orientable 3-manifold $M$ properly dominates at most finitely many closed orientable geometric 3-manifolds with non-empty boundary.
\end{corollary}

Then the following  step is to study Question \ref{lf} when the targets are Haken manifolds (a Haken manifold is not geometric in general but it admits a decomposition into geometric 3-manifolds). This is the purpose of Theorem \ref{dom}.

We end this section by giving an interpretation of Theorem \ref{poset} for the subclass ${\mathcal G}_0$ of ${\mathcal G}$ which consists of graph manifolds.  The purpose of this remark is to study the \emph{local finiteness} of the poset $({\mathcal G}_0,\geq_1)$, up to homotopy equivalence. Recall that a poset $({\mathcal P},\geq)$ is \emph{locally finite} if for any $x,y$ in ${\mathcal P}$ with $x\leq y$ the interval $[x,y]=\{z\in{\mathcal P}, x\leq z\leq y\}$ is finite (many results on posets require this condition). Then Theorem \ref{poset} implies the following

\begin{corollary} The  poset of closed orientable graph manifolds partially ordered, up to homotopy equivalence, by degree one maps is locally finite.
\end{corollary}

\section{Notation and terminology}
\subsection{The degree of a map}
 Let $f\co M\to N$ be a map between orientable compact connected $n$-manifolds. We say that $f$ is proper if $f^{-1}(\b N)=\b M$. Suppose $f$ is proper. Then $f$ induces homomorphisms $f_{\ast}\co\pi_1M\to\pi_1N$, $f_{\sharp}\co H_{\ast}(M,\b M)\to H_{\ast}(N,\b N)$, $f^{\sharp}\co H^{\ast}(N;{\R})\to H^{\ast}(M;{\R})$. The degree of $f$, ${\rm deg}(f)$, is given by the equation $f_{\sharp}([M])={\rm deg}(f)[N]$, where $[M]\in H_n(M,\b M;{\Z})$, $[N]\in H_n(N,\b N;{\Z})$ are the chosen fundamental classes of $M$ and $N$. On the other hand the Gromov simplicial volume ${\rm Vol}(M)$ of the pair $(M,\b M)$ is the infimum of the $l^1$-norms $\sum_{j=1}^{k}|\lambda_j|$ of all cycles $z=\sum_{j=1}^k\lambda_j\sigma_j$, with $\sigma_j\co\Delta^n\to M$ singular $n$-simplexes of $M$, $\lambda_j\in{\R}$, representing the fundamental class  $[M]\in H_n(M,\b M;{\Z})$ (see \cite[Sect. 1.1]{G}). We recall the following well known and useful result on nonzero degree maps.
\begin{proposition}\label{obvious} Suppose $f\co M\to N$ is a proper nonzero degree map between compact orientable $3$-manifolds. Then the following properties hold:

(i) the index of $f_{\ast}(\pi_1M)$ in $\pi_1N$ divides ${\rm deg}(f)$,

(ii) the induced homomorphism $f_{\sharp}\co H_{\ast}(M,\b M;{\R})\to H_{\ast}(N,\b N;{\R})$ is surjective  and by duality $f^{\sharp}\co H^{\ast}(N;{\R})\to H^{\ast}(M;{\R})$ is a monomorphism,

(iii) ${\rm Vol}(M)\geq{\rm deg}(f){\rm Vol}(N)$.
\end{proposition}
\begin{proof}[Sketch of proof]
Point (i) comes directly from a covering space argument as in the proof of Lemma 15.12 in \cite{He}. Point (ii) comes from the Poincar\'e Duality combined with the naturality of cap products. Point (iii) can be obtained directly using the definition of Gromov simplicial volume combined with the definition of the degree of a map given in paragraph 2.1.
\end{proof}
\subsection{Haken manifolds and sewing involutions}
An orientable compact irreducible 3-manifold is called a \emph{Haken manifold} if it contains an orientable proper incompressible surface. Given a closed Haken manifold $N$ we denote by ${\cal T}_N$ the Jaco-Shalen-Johannson family of canonical tori of $N$ and by ${\cal H}(N)$ (resp. ${\cal S}(N)$) the disjoint union of the hyperbolic (resp. Seifert) components of $N\setminus {\cal T}_N\times[-1,1]$ so that $N\setminus{\cal T}_N\times[-1,1]={\cal H}(N)\cup{\mathcal S}(N)$, where ${\cal T}_N\times[-1,1]$ is identified with a regular neighborhood of ${\cal T}_N$ in such a way that ${\cal T}_N\simeq{\cal T}_N\times\{0\}$  (see \cite{JS}, \cite{J} and \cite{T} for the statement and the proof  of this decomposition). On the other hand, we denote by $\Sigma(N)$ the disjoint union of ${\mathcal S}(N)$ with the components of ${\cal T}_N\times[-1,1]$.

Let $N$ be a Haken manifold. Consider the 3-manifold
$N^{\ast}$ obtained after splitting $N$ along ${\mathcal T}_N$. There is an involution $s\co\b N^{\ast}\to\b N^{\ast}$
defined as follows. Let $r\co N^{\ast}\to N$ be the canonical identification map. For any component $T$ of $\b N^{\ast}$ we
denote by $T'$ the unique component of $\b N^{\ast}$ distinct of $T$ such that $r(T')=r(T)$. Let $s_T\co T\to T'$ be the
unique homeomorphism such that $(r|T')\circ s_T=r|T$. Define $s$ by setting $s|T=s_T$ for any $T\in\b N^{\ast}$. The map $s$ will be termed the \emph{sewing involution} for $N$. 

Consider now two Haken manifolds $N_1$ and $N_2$ with sewing involutions $s_1$ and $s_2$. We say that the two ordered
pairs $(N_1^{\ast},s_1), (N_2^{\ast},s_2)$ are \emph{equivalent}  if there is a homeomorphism $\eta\co N_1^{\ast}\to N_2^{\ast}$ such that $\eta\circ s_1$ and
$s_2\circ\eta$ are isotopic. Using this notation then two Haken manifolds $N_1$ and $N_2$ are homeomorphic if and only if the two ordered pairs   $(N_1^{\ast},s_1)$ and $
(N_2^{\ast},s_2)$ are equivalent. On the other hand we will say, for convenience, that two Haken manifolds $N_1$ and $N_2$ are  \emph{weakly equivalent} if there is a homeomorphism $\eta\co N_1^{\ast}\to N_2^{\ast}$. 

\subsection{Haken manifolds, graph manifolds and simplicial volume}
Recall that it follows from \cite{Th6}  that if $H$ is a complete finite volume hyperbolic manifold then $${\rm Vol}(H)=\frac{{\rm Vol}_{int}(H)}{v_3}$$ where ${\rm Vol}_{int}(H)$ is the volume associated to the complete hyperbolic metric in ${\rm int}(H)$ and $v_3$ is a constant which depends only on the dimension. On the other hand it follows from \cite{G} that ${\rm Vol}(S)=0$ when $S$ is a Seifert fibered space. Then using the Cutting off Theorem of M. Gromov (\cite{G}) we get
$${\rm Vol}(N)=\sum_{H\in{\cal H}(N)}{\rm Vol}(H)$$
A 3-manifold $G$ is termed a \emph{graph manifold} if there is a collection ${\mathcal T}$ of disjoint embedded tori in $G$ such that each component of $G\setminus{\mathcal T}$ is Seifert. 
Note that the Gromov simplicial volume gives a characterization of graph manifolds in the following way:
\begin{theorem}[\cite{Sinv}]\label{inv}
A closed orientable 3-manifold $N$ is a graph manifold if and only if $N$ is an element of ${\mathcal G}$ with zero Gromov simplicial volume.
\end{theorem}

We end this section with the following convenient definition. Given a closed Haken manifold $N$, a zero codimensional submanifold $G$ of $N$ which is the union of some geometric (resp. Seifert) components of $N$ will be termed a \emph{ canonical} (resp.\emph{ graph) submanifold  of} $N$.

\section{Main steps of the proof of Theorem \ref{dom} and statement of the intermediate results}
Let $M$ be a closed orientable 3-manifold and let $N$ be a closed Haken manifold 1-dominated by $M$. First note that we may assume, throughout the proof of Theorem \ref{dom},  that the target  satisfies the following condition:

(I) $N$ is a closed non-geometric Haken manifold.

This condition comes from Theorem \ref{known}. On the other hand the constant $c\in(0,1)$ of Theorem \ref{dom} is given by a result of T. Soma in   \cite[Theorem 1]{Sdom} which implies the following
 
 \begin{theorem}[\cite{Sdom}]\label{constant}
 Let $M$ be a closed orientable 3-manifold. There exists a constant $c\in(0,1)$, which depends only on $M$, satisfying the following property. If $f\co M\to N$ denotes a nonzero degree map to a   closed Haken manifold $N$   whose Gromov simplicial volume satisfies ${\rm Vol}(N)\geq(1-c){\rm Vol}(M)$ then ${\rm Vol}(M)={\rm deg}(f){\rm Vol}(N)$.
 \end{theorem}   
This, in order  to state Theorem \ref{dom} we will prove the following general result on non-geometric closed Haken manifolds.
\begin{proposition}\label{A}
Let $M$ be a closed orientable 3-manifold and let $d$ be striclty positive integer. Then there are at most finitely many closed non-geometric Haken manifolds $N$ such that there exists a  degree-$d$ map $f\co M\to N$ satisfying ${\rm Vol}(M)={\rm deg}(f){\rm Vol}(N)$.
\end{proposition}
The proof of Proposition \ref{A} contains two steps. In the first one, we show that there are at most finitely many homeomorphism classes for $N^{\ast}$ (when $N$ runs over the targets manifolds)  and in the second one, we prove that there are at most finitely many equivalence classes of pairs $(N^{\ast},s)$ where $s$ is the sewing map which produces the target $N$ from its geometric decomposition $N^{\ast}$. We give now the key results of this two steps.
\subsection{First step: Control of the geometric decomposition of the targets}
According to the paragraph above, the purpose of this step is to prove the following result:
\begin{proposition}\label{B}
Let $M$ be a closed orientable 3-manifold and let $d$ be a strictly positive integer. Then there are at most finitely many classes of weakly equivalent   non-geometric  closed Haken manifold $N$  such that there exists a  degree-$d$ map $f\co M\to N$ satisfying ${\rm Vol}(M)={\rm deg}(f){\rm Vol}(N)$.
\end{proposition}
The proof of Proposition \ref{B} depends on the following key result which says that a nonzero degree map $f$ into a Haken manifold $N$ has a kind of \emph{canonical standard form} with respect to the geometric decompostion of $N$.  
\begin{lemma}[Standard Form]\label{haken}
Any closed orientable 3-manifold $M$ admits a finite set  ${\mathcal H}=\{M_1,...,M_k\}$ of closed Haken manifolds satisfying the following property. For any nonzero degree map $g\co M\to N$ into a closed non-geometric Haken manifold $N$  containing no embedded Klein bottles and satisfying ${\rm Vol}(M)={\rm deg}(g){\rm Vol}(N)$ there exists at least one element $M_i$ in ${\mathcal H}$ and a nonzero degree map $f\co M_i\to N$ such that:

(i) ${\rm Vol}(M_i)={\rm deg}(f){\rm Vol}(N)$, and

(ii) $f$ induces a finite covering between ${\mathcal H}(M_i)$ and ${\mathcal H}(N)$, and

(iii) for any geometric component $Q$ in $N^{\ast}$ the preimage $f^{-1}(Q)$ is a canonical submanifold of $M$. 
\end{lemma}
\begin{remark}\label{debile}
It will follow from the proof of Lemma \ref{haken}  that if $Q$ is a Seifert piece of $N$ then  $f^{-1}(Q)$ is a graph submanifold of $M_i$ and if $Q$ is a hyperbolic piece then each geometric component of $f^{-1}(Q)$ is a hyperbolic piece of $M_i$. 
\end{remark}
Recall that in \cite[Key Lemma]{Sdom}, T. Soma proves the following result for complete finite volume hyperbolic 3-manifolds without any condition on the Gromov simplicial volume:

\begin{lemma}[T. Soma]\label{sos} Any closed orientable  3-manifold $M$ admits a finite set ${\cal F}=\{F_1,...,F_n\}$ of 3-manifolds such that for any closed Haken manifold $N$ dominated by $M$ then any component $H$ of ${\mathcal H}(N)$ is properly dominated by at least  one element  $F_i$ of ${\cal F}$.
\end{lemma}
Since a closed Haken manifold contains at most finitely many canonical submanifolds then point (iii) of  Lemma \ref{haken}  gives a version of Lemma \ref{sos} for Seifert fibered manifolds with an additional condition on the Gromov simplicial volume. 
First of all, note that in the proof of Lemma \ref{haken} as well as in the proof of Lemma \ref{sos}, it can be shown that there are no loss of generality assuming that $M$ is a closed Haken manifold. With this assumption, recall that the proof of Soma of Lemma \ref{sos} uses the geometry of the  hyperbolic space and in particular the isotropy of hyperbolic geometry is crucial for ``locally hyperbolizing'' certain simplicial subcomplexes of $M$. This method  can not be adapted in the Seifert case since the geometry is not isotropic (indeed there is an invariant direction corresponding to the Seifert fibration).

In the proof of Lemma \ref{haken} the condition on the Gromov simplicial volume is essential. More precisely the proof of Lemma \ref{haken} is based on the observation that when ${\rm Vol}(M)={\rm deg}(f){\rm Vol}(N)$ then we can ``control''  the ``essential part'' of $f^{-1}({\cal T}_N)$. Actually one can show, up to homotopy, that this essential part is a subfamily of ${\cal T}_M$ which is crucial in our proof since this ensures that the genus of the essential components of $f^{-1}({\cal T}_N)$ is bounded independently of $N$. This control can not be accomplished when ${\rm Vol}(M)>>{\rm  deg}(f){\rm Vol}(N)$. Indeed, consider for example a degree one map from a closed hyperbolic 3-manifold $M$ to a graph manifold $N$ (this kind of example can be built by taking a hyperbolic nul-homotopic knot $k$ in a graph manifold $N$  and by gluing a solid torus along $\partial(N\setminus k)$ in such a way that the resulting manifold $M$ is  hyperbolic, then the degree of the  canonical decomposition map $f:M\to N$ is one, see \cite{BW} for details on this construction).  In this case one can clearly not control the genus of the components of  $f^{-1}({\cal T}_N)$.

The family ${\mathcal H}$ of Haken manifolds in Lemma \ref{haken} comes from a finite family of canonical  submanifolds ${\mathcal A}$ of $M$ after some Dehn fillings. Note that to get a family $\hat{\mathcal A}$ of Haken manifolds whose elements satisfies conditions (i), (ii) and (iii) one can use a construction of Rong in \cite{Ro2}. But this construction does not guarantee the finiteness of the family $\hat{\mathcal A}$ (actually the construction of Rong does not allow to control the slopes of the Dehn fillings performed along the components of ${\mathcal A}$ to obtain $\hat{\mathcal A}$). Thus we have to modify this construction to avoid this problem. To this purpose we will define and  construct the \emph{maximal essential part} of $M$ (see Section 5.3).

\subsection{Second step: Control  of  the sewing involutions of the targets}
In this step we complete the proof of Proposition \ref{A}. Thus the key result of this section states as follows.
\begin{proposition}\label{C}
Let $M$ be a closed orientable 3-manifold and let $d$ be a strictly positive integer. Let $N_i$ be a sequence of weakly equivalent non-geometric closed Haken manifolds   such that there exists a degree-$d$ map $g_i\co M\to N_i$ satisfying ${\rm Vol}(M)={\rm deg}(g_i){\rm Vol}(N_i)$. For each $i\in{\N}$,  we denote by $s_i\co\b N_i^{\ast}\to\b N_i^{\ast}$  the sewing involution corresponding to $N_i$. Then the sequence $\{(N_i^{\ast},s_i), i\in{\N}\}$ is finite, up to equivalence of pairs.  
\end{proposition}
Throughout the proof of Proposition \ref{C} we will use the collection of closed Haken manifolds ${\mathcal H}$ given by Lemma \ref{haken}. Points (i), (ii) and (iii) say that the elements of ${\mathcal H}$ dominate the manifolds $N_i$'s in a convenient way. Roughly speaking, the core of the proof of Proposition \ref{C} is to show that the sewing involution associated to each Haken manifold of ${\mathcal H}$ does fix  the sewing involution $s_i$ which produces $N_i$ from $N_i^{\ast}$. Note that in this step the condition on the Gromov simplicial volume is still crucial in our proof.

\subsection{Organization of the paper} This paper is organized as follows.  Section 4 is devoted to the statement of a \emph{mapping result} for maps from  Seifert fibered spaces to Haken manifolds. This result  has only a technical interest and will be used in Sections 5 and 6.  Section 5 is devoted to the proof of Proposition \ref{B}  and in   Section 6 we prove  Proposition \ref{C} to complete the proof of   Proposition \ref{A}.  Section 7 is devoted to the  proof of   Theorem  \ref{poset} which is a consequence of Theorems \ref{dom} and \ref{known}.

\section{On the Characteristic Pair Theorem of W. Jaco and P. Shalen} 

We start by recalling a main consequence of the Characteristic Pair Theorem of W. Jaco and P. Shalen (see \cite[Chapter V]{JS}) which allows to control a nondegenerate map from a Seifert fibered space into a  Haken manifold.  We first give the definition of degenerate maps in the sense of W. Jaco and P. Shalen.
\begin{definition}\label{deg}
Let $(S,F)$ be a connected Seifert pair, and let $(N,T)$ be a connected 3-manifold pair. A map $f\co(S,F)\to(N,T)$ is said to be \emph{degenerate} if either

(0) the map $f$ is inessential as a map of pairs, or

(1) the group ${\rm Im}(f_{\ast}\co\pi_1S\to\pi_1N)=\{1\}$, or

(2) the group ${\rm Im}(f_{\ast}\co\pi_1S\to\pi_1N)$ is cyclic and $F=\emptyset$, or

(3) the map $f|\gamma$ is homotopic in $N$ to a constant map for some fiber $\gamma$ of $(S,F)$.    
\end{definition}
Then the Characteristic Pair Theorem of Jaco and Shalen implies the following result.
\begin{theorem}\label{ndeg}[Jaco, Shalen]
If $f$ is a nondegenerate map of a Seifert pair $(S,\emptyset)$ into a Haken manifold pair $(M,\emptyset)$, then there exists a map $f_1$ of $S$ into $M$, homotopic to $f$, such that $f_1(S)\subset{\rm int}(\Sigma(M))$.
\end{theorem}
The purpose of this section is to give a kind of \emph{mapping lemma} for a certain class of degenerate maps. More precisely we show here the following result which will be used in the proof of Theorem \ref{dom}.
 \begin{lemma}\label{mapping}
 Let $f\co M\to N$ be a  map between closed Haken manifolds and suppose that $N$ is non-geometric and contains no embedded Klein bottles.  Let $S$ and $S'$ be two components of ${\mathcal S}(M)$ which are adjacent in $M$ along a subfamily ${\mathcal T}$ of ${\mathcal T}_M$. Assume that $S$ and $S'$ satisfy the following hypothesis:
 
 (i) $f(S')\subset{\rm int}(B')$, where $B'$ is a component of $\Sigma(N)$, and
 
 (ii) $f_{\ast}(t_S)\not=1$, where $t_S$ denotes the homotopy class of the regular fiber of $S$.
 
 Then there exists a component $B$ of $\Sigma(N)$, with regular fiber $h$, and a homotopy $(f_t)_{0\leq t\leq 1}$ which is constant outside of a regular neighborhood of $S$ such that $f_0=f$ and $f_1(S)\subset{\rm int}(B)$. Moreover if $(f_1)_{\ast}(t_S)$ is not conjugate to a non-trivial power of   $h$  then one can choose $B=B'$ and  thus $f_1(S\cup_{\mathcal T}S')\subset{\rm int}(B')$.
 \end{lemma}
 \begin{figure}[h]
\centerline{\input{mapping.pstex_t}}
\end{figure}
 \begin{proof}
 Let $T$ be a canonical torus of $M$ such that $T\in\b S\cap\b S'$ and denote by $t_S$ the regular fiber of $S$ represented in $T$. It follows from the hypothesis of the lemma that there exists a Seifert piece $B'$ of $\Sigma(N)$ such that  $f(S')\subset B'$ and thus $f_{\ast}(t_S)\in\pi_1B'\setminus\{1\}$. Fix a base point $x$ in $T$, in such a way that the groups  $\pi_1S$ and $\pi_1S'$ are always considered with base point $x$ and denote by $y=f(x)$ a base point in $B'$.

\emph{Case 1.} If $f_{\ast}(\pi_1S)$ is nonabelian, since $f_{\ast}(t_S)\not=\{1\}$, then $f|S\co S\to N$ is a nondegenerate map. Hence   the Characteristic Pair Theorem implies that there exists $B\in\Sigma(N)$ such that $f(S)\subset{\rm int}(B)$. Moreover since $f_{\ast}(\pi_1S)$ is nonabelian then $f_{\ast}(t_S)$ has nonabelian centralizer and \cite[Addendum to Theorem VI.I.6]{JS} implies that $f_{\ast}(t_S)\in\l h\r$, where $h$ denotes the regular fiber of $B$. This proves the lemma when $f_{\ast}(\pi_1S)$ is nonabelian.
 
 Assume that $f_{\ast}(\pi_1S)$ is abelian. Since $\pi_1N$ is torsion free, and since $N$ is an aspherical 3-manifold then the subgroup $f_{\ast}(\pi_1S)$ of $\pi_1N$ must have cohomological dimension at most 3 and thus it is isomorphic to either  ${\Z}$ or ${\Z}\times{\Z}$ or ${\Z}\times{\Z}\times{\Z}$. The case ${\Z}\times{\Z}\times{\Z}$ is excluded since $N$ is a non-geometric closed Haken manifold.
 
\emph{Case 2.} Thus assume first that $f_{\ast}(\pi_1S)\simeq{\Z}\times{\Z}$. In this case $f|S\co S\to N$ is still  a nondegenerate map and the Characteristic Pair Theorem implies that there exists component $B\in\Sigma(N)$, with regular fiber $h$, adjacent to $B'$ in $N$ such that $f(S)\subset{\rm int}(B)$, after a homotopy on $f$.  Suppose that  $f_{\ast}(t_S)\not\in\l h\r$. Thus by \cite[Addendum to Theorem VI.I.6]{JS} we know that the centralizer ${\mathcal Z}_{|\pi_1(B,y)}(f_{\ast}(t_S))$ of $f_{\ast}(t_S)$ in $\pi_1(B,y)$ is necessarily abelian. Let $c$ be an element of $\pi_1S$. Then $f_{\ast}(c)\in{\mathcal Z}_{|\pi_1(B,y)}(f_{\ast}(t_S))$. Denote by $h'$ the regular fiber of $B'$ represented in a component of $B\cap B'$ in such a way that 
$$h'\in\pi_1(B,y)\cap\pi_1(B',y)$$ 
Since $f_{\ast}(t_S)\in\pi_1(B',y)\cap\pi_1(B,y)$ (recall that $t_S\in\pi_1(S,x)\cap\pi_1(S',x)$) then $h'$ commutes with $f_{\ast}(t_S)$ and since $h'\in\pi_1(B',y)\cap\pi_1(B,y)$ then  $h'\in{\mathcal Z}_{|\pi_1(B,y)}(f_{\ast}(t_S))$.  Thus, since ${\mathcal Z}_{|\pi_1(B,y)}(f_{\ast}(t_S))$ is abelian  this implies that $f_{\ast}(c)\in Z(h')$. Since $c$ is an arbitrary element in $\pi_1S$ then $f_{\ast}(\pi_1S)\subset Z(h')$. This implies that $f_{\ast}(\pi_1S)$ is conjugate to a subgroup of $\pi_1(B',y)$. Then after a homotopy on $f$ we may assume that $f(S)\subset{\rm int}(B')$. This prove the lemma when $f_{\ast}(\pi_1S)\simeq{\Z}\times{\Z}$.

\emph{Case 3.}  Assume now that  $f_{\ast}(\pi_1S)\simeq{\Z}$. Then there exists an element $c\in\pi_1S$ such that $f_{\ast}(\pi_1S)=\l f_{\ast}(c)\r$ and in particular there exists $n\in{\Z}^{\ast}$ such that $f_{\ast}(t_S)=(f_{\ast}(c))^n$.  In the following  $[a,b]$ denotes the commutator of $a$ and $b$. Since in this case the Characteristic Pair Theorem does not apply, since $f|S\co S\to N$ is a degenerate map, we  first prove that there exists $B\in{\mathcal S}(N)$ such that $f(S)\subset{\rm int}(B)$, after a homotopy on $f$.
 
  \emph{Subcase 3.1.} Assume that $[f_{\ast}(c),h']=1$. In this case $f_{\ast}(c)$, and hence $f_{\ast}(\pi_1S)$, is in the centralizer of $h'$ and thus  one can deform $f$ on a regular neighborhood of $S$ such that $f(S)\subset{\rm int}(B')$.  
 
 \emph{Subcase 3.2.} Assume that $[f_{\ast}(c),h']\not=1$. Since $f_{\ast}(c)$ and $h'$ are in the centralizer  $Z(f_{\ast}(t_S))$ of $f_{\ast}(t_S)$ then the group $Z(f_{\ast}(t_S))$ is non-abelian. Then by \cite[Addendum to Theorem VI.I.6]{JS} we know that $f_{\ast}(t_S)$ is conjugate to a power of the regular fiber $h$ of a Seifert piece $B$ of ${\mathcal S}(N)$. Thus one can deform $f$ on a regular neighborhood of $S$ such that $f(S)\subset{\rm int}(B)$. Note that since a power of $f_{\ast}(c)$ lies in $\l h\r$ then by \cite[Lemma II.4.2]{JS}, $f_{\ast}(c)=c_i^{\alpha_i}$, where $c_i$ denotes the homotopy class of an exceptional fiber in $B$ and ${\alpha_i}\in{\Z}^{\ast}$. 
 
 To complete the proof of the lemma in Case 3 it is sufficient to apply the same argument as in case 2. 

 \end{proof}

\section{Control of the geometric pieces of the targets}
This section is devoted to the proof of Proposition \ref{B}. To this purpose we first give a proof of Lemma \ref{haken}. Let $M$ be a closed orientable 3-manifold and let  $f\co M\to N$ be a nonzero degree map into a closed non-geometric Haken manifold which contains no embedded Klein Bottles such that ${\rm Vol}(M)={\rm deg}(f){\rm Vol}(N)$. First we claim that to prove Lemma \ref{haken} there is no loss of generality assuming that $M$ is a closed Haken manifold. Indeed, consider the Milnor decomposition of $M$ into prime manifolds $M=M_1\sharp...\sharp M_k$ (see \cite{M}). Since $\pi_2(N)$ is trivial, there exists, for each $i\in\{1,...,k\}$ a map $f_i\co M_i\to N$ such that ${\rm deg}(f_1)+...+{\rm deg}(f_k)={\rm deg}(f)$. Note that when ${\rm deg}(f_i)\not=0$ then $M_i$ is necessarily a closed Haken manifold. On the other hand if ${\rm Vol}(M)={\rm deg}(f){\rm Vol}(N)$  then the Cutting of Theorem of M. Gromov, \cite{G}, implies that there exists $i\in\{1,...,k\}$ such that $f_i\co M_i\to N$ has nonzero degree and satisfies ${\rm Vol}(M_i)={\rm deg}(f_i){\rm Vol}(N)$. Then from now one we assume that $M$ is a closed Haken manifold.
\subsection{A convenient alternative  to  Lemma \ref{haken}.}

\subsubsection{Sections of Seifert fibered spaces}
Let $S$ be an orientable ${\Hi}^2\times{\R}$-Seifert fibered space with non-empty boundary and orientable basis $B$. Then the Seifert fibration of $S$ is unique and we denote by $\eta\co S\to B$ the canonical projection map. If $S$ has exceptional fibers $C_1,...,C_r$, let $D_1,...,D_r$ be pairwise disjoint 2-cells neighborhood of $\eta(C_1),...,\eta(C_r)$ in ${\rm int}(B)$. Let $B'=B\setminus\cup_i{\rm int}(D_i)$ and $S'=\eta^{-1}(B')$. Then $\eta|S'\co S'\to B'$ is the orientable circle bundle over $B'$ and since $B'$ is orientable then $S'=B'\times{\S}^1$. Choose a cross section $s_0\co B'\to S'$ of the circle bundle. We may choose standards generators of $\b S'$, with respect to this choice of a cross section, in the following way. Denote $\b S'=\b S\cup U_1\cup ...\cup U_r$ where $U_j=\b\eta^{-1}(D_j)$. Then for each component $U_j$ (resp. $T_i$ of $\b S$) we choose generators $t$, $q_j$ (resp. $t$, $\delta(S,T_i)$) where $t$ is represented by a regular fiber and $q_j$ (resp. $\delta(S,T_i)$) is the boundary curve of the cross section $s_0$ in $U_j$ (resp. in $T_i$). In the following the curve $\delta(S,T_i)$ will be termed  a section of $T_i$ (with respect to the Seifert fibration of $S$). Notice that if we replace the section $s_0$ by an other one $s\co  B'\to S'$ then the section $\delta(S,T_i)$ of $T_i$ is replaced by $\delta(S,T_i)t^m$, $m\in{\Z}$.

\subsubsection{Dehn fillings} 
Let $Q$ be a compact oriented three manifold whose boundary is made of tori $T_1,...,T_k$. For each $i=1,...,k$ we fix generators $l_i,m_i$ of $\pi_1T_i$. Let ${\cal P}_{\ast}$ be the subset of ${\bf S}^2={\C}\cup\{\infty\}$ defined by $${\cal P}_{\ast}=\{(p,q)\in{\Z}\times{\Z}, {\rm gcd}(p,q)=1\}\cup\{\infty\}$$ where ${\rm gcd}(p,q)$ denotes the greatest common divisor of $p$ and $q$. We will denote by $Q_{d_1,...,d_k}$ the 3-manifold obtained from $Q$ by gluing to each $T_i$, $i=1,...,k$, a solid torus $S^1\times D^2$ identifying a meridian $m=\{z_0\}\times\b D^2$ with $p_il_i+q_im_i$ when $d_i=(p_i,q_i)\in{\mathcal P}_{\ast}\setminus\{\infty\}$. When $d_i=\infty$ the torus $T_i$ is cut out. On the other hand recall that the manifolds obtained in this way depend, up to diffeomorphism, only on the pair of integers $(p_i,q_i)$ with ${\rm gcd}(p_i,q_i)=1$. Let $M$ be  closed Haken manifold. From now on we adopt the following convention.

 For each $T$ in $\b{\cal S}(M)$ we fix a Seifert fibered space $S$ adjacent to $T$ and a basis $(h_T,\delta_T)$ of $\pi_1(T)$ where $h_T$ corresponds to the generic fiber $h(S)$ of $S$ and $\delta_T$ is a section $\delta(S,T)$ of $T$ with respect to the Seifert fibration of   $S$ as defined in Paragraph 4.1.1.  If $S$ is adjacent to a Seifert fibered space  $S'$ along $T$ we denote by $(h(S'),\delta(S',T))$ an other basis for $\pi_1T$ with respect to $S'$ in the same way as for $S$. We denote by $d_T=(a_T,b_T)$  the element of ${\cal P}_{\ast}$ such that $h(S')=a_Th_T+b_T\delta_T$.
 Note that $b_T\not=0$ by the minimality property  of ${\cal T}_M$. Denote by ${\cal P}_{\ast}^0$ the finite subset of ${\cal P}_{\ast}$ defined  by
$${\cal P}_{\ast}^0=\{(a_T,b_T), T\in\b{\cal S}(M)\setminus\b{\cal S}(M)\cap\b{\cal H}(M),  (1,0)\}$$        Then to prove Lemma \ref{haken} it is sufficient to  state the following result.

\begin{lemma}\label{splitting} Let $M$ be a closed Haken manifold and let $N$ be a  closed non-geometric Haken manifold  that contains no embedded Klein bottles. If    $f\co M\to N$ denotes a nonzero degree map satisfying ${\rm Vol}(M)={\rm deg}(f){\rm Vol}(N)$ then there exists a canonical submanifold $G_N$ of $M$ whose boundary is made of some components of  $\b{\cal S}(M)\setminus\b{\cal S}(M)\cap\b{\cal H}(M)$ and such that if $T_1,...,T_k$ denotes the components of $\b G_N$ then there exists $d_1,....,d_k$  in  ${\cal P}_{\ast}^0$  satisfying the following properties:

(a)  $(G_N)_{d_1,...,d_k}$ is a closed Haken manifold, and

(b) there exists a nonzero degree map $g\co (G_N)_{d_1,...,d_k}\to N$ satisfying points (i), (ii) and (iii) of Lemma \ref{haken}.
\end{lemma}

\subsection{Non-zero degree maps preserving the Seifert part of the domain.}
In this section we prove that Lemma \ref{splitting} is true for non-zero degree maps $f:M\to N$ such that ${\rm Vol}(M)={\rm deg}(f){\rm Vol}(N)$ and satisfying  $f({\cal S}(M))\subset{\rm int}({\cal S}(N))$.
\begin{lemma}\label{splittingregular} Let $f:M\to N$ be a nonzero degree map between non-geometric Haken manifolds such that ${\rm Vol}(M)={\rm deg}(f){\rm Vol}(N)$. If $f({\cal S}(M))\subset{\rm int}({\cal S}(N))$ then there exists a map homotopic to $f$ which satisfies the conclusion of Lemma \ref{haken}.
\end{lemma}

 \begin{proof} First of all note that using the construction of T. Soma in \cite{Srigide} one can modify $f$ by a homotopy fixing $f|{\cal S}(M)$ in such a way that $f({\cal H}(M),\b{\cal H}(M))\subset({\cal H}(N),\b{\cal H}(N))$ and $f|{\cal H}(M):{\cal H}(M)\to{\cal H}(N)$ is a ${\rm deg}(f)$-fold covering.

Let $T\in{\cal T}_N$. Using standard cut and paste arguments and
the fact that $\b{\cal S}(M)$ and $\b{\cal H}(M)$ are incompressible
we can modify $f$ by a homotopy fixing $f|{\cal S}(M)\cup{\cal
H}(M)$, so that $f^{-1}(T)$ is a collection of 2-sided
incompressible surfaces in $M$. Since $f^{-1}(T)\subset M\setminus({\cal
S}(M)\cup{\cal H}(M))$ it must be a union of parallel copies of
some tori in ${\cal T}_M\times(-1,1)$. We can arrange $f$ in its
homotopy class so that for any $U\in{\cal T}_M$, a regular neighborhood $U\times[-1,1]$ of $U$  contains at most one
component of $f^{-1}(T)$. Indeed, suppose that $X$ and $X'$ are
two consecutive components of $f^{-1}(T)\cap(U\times[-1,1])$. Then
$X$ and $X'$ bound a region $Q$ in $U\times[-1,1]$ which is
homeomorphic to ${\S}^1\times{\S}^1\times I$ and there is a
Seifert piece $B$ in ${\cal S}(N)$  so that $f(Q,\b Q)\subset(B,\b
B)$. Then by \cite[Lemma 2.8]{Ro2}, $f|Q$ is homotopic, mod. $\b Q$,
to a map $f_1$ such that $f_1(Q)\subset\b B$, unless
$B\simeq{\S}^1\times{\S}^1\times I$ which is excluded since $N$ is
not a geometric 3-manifold. So we can eliminate
$X$ and $X'$ by pushing $Q$ into $N-B$. After repeating this
operation a finite number of times we may assume that
$f^{-1}(T)\cap(U\times[-1,1])$ has at most one component.

 Note that since $f:M\to N$ is a non zero degree map then $f_{\ast}(\pi_1M)$ has finite index in $\pi_1N$ and thus for any $S$ in ${\cal S}(N)$ there exists at least   one component of ${\cal S}(M)$ which is sent into ${\rm int}(S)$ via $f$. So $f^{-1}(S)$ consists of some components of ${\cal S}(M)$ union some $T\times[-1,1]$ for $T$ in ${\cal T}_M$ (precisely when $f^{-1}({\cal T}_N)\cap(T\times[-1,1])=\emptyset$). So each component of $f^{-1}(S)$ is a canonical graph submanifold of $M$.     This proves Lemma \ref{splittingregular}.
 \end{proof}

 \subsection{Proof of  Lemma \ref{splitting}: the general case.}
 We first realize a kind of {\it factorization} on the map $f$ which is inspired from a construction of Y. Rong in \cite{Ro2}  to have a  reduction to the case of Lemma \ref{splittingregular}. If $S$ is a component of ${\cal S}(M)$ we denote by $t_S$ the homotopy class of the regular fiber in $S$.   Let $B_0$ be the union of all $S$ in ${\cal S}(M)$
such that $f|S$ is degenerate in the sense of Definition \ref{deg}.  If $f|S$ is
 a degenerate map then  either 
 \begin{enumerate}
 \item[Case 1]{ $f_{\ast}(\pi_1 S)=\{1\}$  or,}
 
\item[Case 2]{$f_{\ast}(\pi_1 S)={\Z}$  or,} 

\item[Case 3]{Since $\pi_1(N)$ is torsion free,  $(f|S)_{\ast}:\pi_1 S\to\pi_1
N$ factors through $\pi_1 V$, where $V$ is the base 2-manifold of
the Seifert fibered space $S$.}
\end{enumerate}
 Set $G_0=\o{M-B_0}$. Define a subset of $B_0$ by setting
$${\cal S}_0=\{S\in B_0\setminus(B_0\cap{\mathcal T}_M) \ \mbox{s.t.} \ S \ \mbox{is adjacent to}\ G_0\ \mbox{and}\ f_{\ast}(t_S)\not=1\}$$
and set $B_1=B_0-{\cal S}_0$ and $G_1=\o{M-B_1}$.  We continue this process by setting $${\cal S}_1=\{S\in B_1\setminus(B_1\cap{\mathcal T}_M) \ \mbox{s.t.} \ S \ \mbox{is adjacent to}\ G_1\ \mbox{and}\ f_{\ast}(t_S)\not=1\}$$ to construct an increasing sequence $G_0\subset G_1\subset...\subset G_i\subset G_{i+1}\subset...$ of canonical submanifolds of $M$. We claim that this sequence satisfies the following conditions:
\begin{enumerate}
\item{  the number of connected components $n_{i}$ of $G_{i}$ satisfies  $n_{i+1}\leq n_i$,}

\item{ for any $i$, int$(G_i)$ contains ${\cal H}(M)$ and $f|\b{\mathcal H}(M)\co\b{\mathcal H}(M)\to N$ is a non-degenerate (i.e. $\pi_1$-injective) map,}

\item{ for any $i$ there exists a non-zero degree map $\beta_i:\hat{G_i}\to N$ such that deg$(\beta_i)=$ deg$(f)$, where $\hat{G_i}$ denotes the space obtained from $G_i$ after performing some Dehn fillings along the components of $\b G_i$.}
\end{enumerate}
\begin{figure}[h]
\centerline{\input{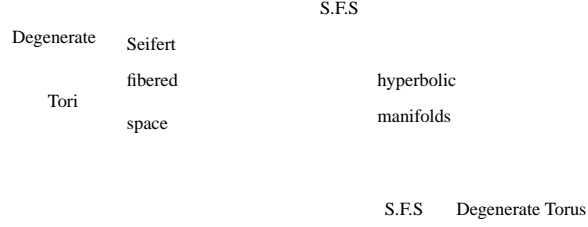}}
\caption{Essential submanifold $G_0$ of $M$ with respect to $f:M\to N$}
\end{figure}
For this reason  $G_i$ is called an \emph{essential part} of $M$ with respect to $f$.
We define an integer $n_0$ by setting:  $$n_0=\min\{n\geq 0\ \mbox{such that}\ G_n=G_{n+1}\}$$  We  prove point    (3) for $G_{n_0}$ which will be termed a \emph{maximal essential part} of $M$. The proof for the $G_i\subset G_{n_0}$ works in the same way and the proof of points (1) and (2) follows directly from the construction and from the equality ${\rm Vol}(M)={\rm deg}(f){\rm Vol}(N)$ (see \cite[Lemma 2]{Srigide}).
  Denote by $B_{\Z}$ the subset of $G_{n_0}\setminus(G_{n_0}\cap{\cal T}_M)$ which consists of the Seifert fibered spaces which are degenerate under $f$.  Note that it follows from the construction that for any $S$ in $B_{\Z}$, $f_{\ast}(\pi_1S)$ is necessarily infinite cyclic. Set $B=B_0-B_{\Z}$. We have $G_{n_0}=\o{M-B}$.

 Let $Q$ be a geometric piece in $G_{n_0}$ such that $\b Q\cap\b G_{n_0}\not=\emptyset$. Then it follows from the construction   that  $Q$ is a Seifert fibered space  and it is adjacent along each component of $\b Q\cap\b G_{n_0}$ to a degenerate Seifert piece in $M$ whose fibers are sent trivially in $\pi_1N$.
For any $S$ in $B$,   define a group $\pi_S$ to be
  \begin{enumerate}
\item[Case 1]{  
$\{1\}$ or,}
\item[Case 2]{ ${\Z}$  or,} 
\item[Case 3]{$\pi_1V$}
\end{enumerate}
 and a
three dimensional space $D_S=K(\pi_S,1)$. Since $D_S$ and $N$ are
both $K(\pi,1)$ there exist maps $\alpha:S\to D_S$ and
$\beta:D_S\to N$ such that $f|S$ is homotopic to
$\beta\circ\alpha$ and satisfying the following convenient conditions: for
each $T\subset\b S$, let $\{\lambda, \mu\}$ be a base of $\pi_1T$
with $\alpha_{\ast}(\lambda)=1$. Note that it follows from the construction that for any $T$ in $\b B=\b G_{n_0}$ then  $\lambda=h_S(T)$, where $S$ is the Seifert fibered manifold of $B$ containing $T$ in its boundary and where $h_S(T)$ denotes the regular fiber of $S$ represented in $T$. Parametrize $T$ by
$T={\S}^1\times{\S}^1$ with $[{\S}^1\times\ast]=\lambda$ and
$[\ast\times{\S}^1]=\mu$. Then $\alpha(x,y)=\alpha_1(y)$ for some
embedding $\alpha_1:{\S}^1\to D_S$. Denote the knot
$\alpha_1({\S}^1)$ by $l_T$. We may also assume that $l_{T_1}\cap
l_{T_2}=\emptyset$ for different components $T_1$ and $T_2$ of $\b
S$. We extend the homotopy on $f|S$ over $M$, we replace $f$ by
the new map and we do this for each component $S$ of $B$. Set
$D_B=\bigcup_{S\in B}D_S$. Then the following diagram commutes:
$$\xymatrix{
B \ar[r]^{f|B} \ar[d]_{\alpha} & N\\
D_B \ar[ur]_{\beta}}$$
Let $\hat{G}_{n_0}$ be the closed 3-manifold
obtained from $G_{n_0}$ by attaching a solid torus $V_T$ to $G_{n_0}$ along each
component $T$ of $\b G_{n_0}=\b B$ so that the meridian of $V_T$ is
identified with the curve  $\lambda$ defined above. Let $l_T'$ be the core of $V$ which has
the same orientation as $\mu$. Let $X=D_B\cup_{\tau}\hat{G}_{n_0}$ where
$\tau$ identifies each $l_T$ in $D_B$ with $l_T'$ in $\hat{G}_{n_0}$
(preserving orientation). Define
$\alpha:G_{n_0}\to\hat{G}_{n_0}=G_{n_0}\cup_{T\in\b G_{n_0}}V_T$ to be the map such
that $$\alpha|G_{n_0}\setminus\b G_{n_0}:G_{n_0}\setminus\b G_{n_0}\to\hat{G}_{n_0}\setminus\cup l_T'$$ is a
homeomorphism and each $T\subset\b G_{n_0}$ is sent onto $l_T'$. Now the
map $\alpha:M=B\cup G_{n_0}\to X$ is a well defined continuous map.
Since $\alpha|G_{n_0}\setminus\b G_{n_0}:G_{n_0}\setminus\b G_{n_0}\to\hat{G}_{n_0}\setminus\cup l_T'$ is a
homeomorphism we can define $\beta|{\hat{G}_{n_0}\setminus\cup
l_T'}=f\circ\alpha_{|\hat{G}_{n_0}\setminus\cup l_T'}^{-1}$. So we get a map $\beta:X\to N$ such that
the following diagram commutes:
$$\xymatrix{
M=B\cup G_{n_0} \ar[r]^{f} \ar[d]_{\alpha} & N\\
X=D_B\cup_{\tau}\hat{G}_{n_0} \ar[ur]_{\beta}}$$
More precisely let $T_1,...,T_l$ be the components of $\b G_{n_0}=\b B$ and let $S_1,...,S_l$ (resp. $B_1,...,B_l$) be the Seifert pieces (not necessarily pairwise distinct) in $G_{n_0}$ (resp. in $B$) such that for each $i=1,...,l$, $B_i$ and $S_i$ are adacent along $T_i$. Denote by 
$$(h(B_i),\delta(B_i,T_i))\ \ {\rm resp.}\ \ (h(S_i),\delta(S_i,T_i))$$
 a system of generators of $\pi_1T_i$ where $h(B_i)$ (resp. $h(S_i)$) denotes the generic fiber of $B_i$ (resp. $S_i$) represented in $T_i$ and $\delta(B_i,T_i))$ (resp. $\delta(S_i,T_i))$) is a section of $T_i$ (with respect to $B_i$, resp. $S_i$) as defined in Section 4.1.1. We know from the construction that $f_{\ast}(h(B_i))=1$. Let $(a_{T_i}(S_i),b_{T_i}(S_i))$ denote the element of ${\cal P}_{\ast}$ such that $$h(B_i)=a_{T_i}(S_i).h(S_i)+b_{T_i}(S_i).\delta(S_i,T_i)$$ If $(h_{T_i},\delta_{T_i})=(h(B_i),\delta(B_i,T_i))$ we set $d_i=(1,0)\in{\cal P}_{\ast}^0$ and if $(h_{T_i},\delta_{T_i})=(h(S_i),\delta(S_i,T_i))$  then $(a_{T_i}(S_i),b_{T_i}(S_i))=(a_{T_i},b_{T_i})\in{\cal P}_{\ast}^0$ and we set $d_i=(a_{T_i},b_{T_i})\in{\cal P}_{\ast}^0$ (see paragraph 4.1.2 for the notations).  Thus we get $\hat{G}_{n_0}=(G_{n_0})_{d_1,...,d_l}$.
\begin{figure}[h]
\centerline{\input{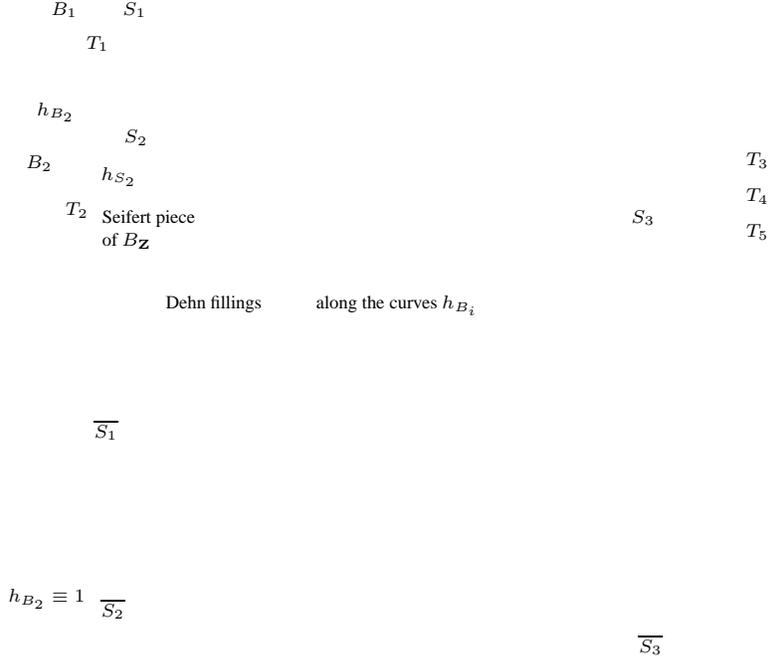}}
\caption{Maximal essential submanifold $G=G_{n_0}$ of $M$}
\end{figure}
Denote by $f_{1}$ the map
$\beta\circ i\co\hat{G}_{n_0}\to N$ where $i\co\hat{G}_{n_0}\to X$ is the
inclusion. Note that since $H_3(D_B)=0$ then a Mayer Vietoris argument shows that $f_1$ is a non-zero degree
map equal to ${\rm deg}(f)$.
\begin{remark}\label{alacon} Let $G^1,...,G^m$ be the components of $\hat{G}_{n_0}$. Up to re-indexing we may assume that there exists $1\leq u\leq m$ such that deg$(f_1|G^i)\not=0$ for $i=1,...,u$ and  deg$(f_1|G^i)=0$ for $i=u+1,...,m$. Set $\hat{G}=G^1\cup...\cup G^u$. There is no loss of generality assuming that $\hat{G}_{n_0}=\hat{G}$.
\end{remark}
 Thus to complete the proof of Lemma \ref{splitting} it remains to check, in view of Lemma \ref{splittingregular}, the following claim.

\begin{claim}\label{prout} The space $\hat{G}_{n_0}$ is a Haken manifold and the map $f_1:\hat{G}_{n_0}\to N$  satisfies ${\rm Vol}(\hat{G}_{n_0})={\rm deg}(f_1){\rm Vol}(N)$. Moreover there exists a map $g:\hat{G}_{n_0}\to N$ homotopic to $f_1$  such that  $g({\cal S}(\hat{G}_{n_0}))\subset{\rm int}({\cal S}(N))$.
\end{claim}
\begin{proof}[Proof of Claim]
Each component $G^i$  of $\hat{G}_{n_0}$ is $\hat{M}_i=(M_i)_{d_1,...,d_{j_i}}$ where $M_i$ is a
union of some hyperbolic pieces and some Seifert fibered pieces of
$M$ connected by some $T\times I$ in ${\cal T}_M\times I$.  Note
that it follows from the construction that for each $i$ the minimal torus
decomposition of $M_i$ gives in an obvious way the minimal torus decomposition of  $G^i=\hat{M_i}$ in the sense that there exists a subfamily ${\cal T}$ of ${\cal T}_M\cap{\rm int}(G_{n_0})$ such that $\alpha({\cal T})={\cal T}_{\hat{G}_{n_0}}$ (recall that $\alpha|G_{n_0}\setminus\b G_{n_0}:G_{n_0}\setminus\b G_{n_0}\to\hat{G}_{n_0}\setminus\cup l_T'$ is a homeomorphism).
We describe precisely the torus decomposition of $\hat{G}_{n_0}$.

{\it Case 1:} Let $S$ be a component of ${\cal S}(M)\cap G_{n_0}$ such that $\b S\cap\b G_{n_0}\not=\emptyset$ and such that $f|S:S\to N$ is a non degenerate map. Set $\hat{S}=\alpha(S)$ in $\hat{G}_{n_0}$. Then $\hat{S}$ admits a Seifert fibration which extends that of $S\subset G_{n_0}$ and $\hat{S}$ is not homeomorphic to a solid torus, $\b\hat{S}$ (if non empty) is incompressible and $f_1|\hat{S}:\hat{S}\to N$ is a non degenerate map. Indeed let $T$ be a component of $\b S\cap\b G_{n_0}$ and let $\lambda$ be the primitive curve of $T$ defined as before. Since $\lambda$ is not a fiber of $S$, by the definition of non-degenerate maps, then the Seifert fibration of $S$ extends to a Seifert fibration of $S\cup_{\lambda=m}V_T$, where $V_T$ denotes a solid torus glued along $T$ by identifying $\lambda$ with the meridian $m$ of $V_T$. Now since $\pi_1(\hat{S})=\pi_1S/\l\lambda\r$ maps onto $\pi_1S/\ker(f_{\ast})\simeq f_{\ast}(\pi_1S)$ which is not cyclic by the definition of non-degenerate maps, then $\hat{S}$ is not a solid torus and so $\b\hat{S}$ is incompressible. Moreover notice that if a torus $T$ connects two non-degenerate Seifert piece $S_1$ and $S_2$ in $M_i$ then $T$ also connects $\alpha(S_1)=\hat{S}_1$ and $\alpha(S_2)=\hat{S}_2$ and the fibers of $\hat{S}_1$ and $\hat{S}_2$ do not match up along $T$ and thus $T\in{\cal T}_{\hat{G}_{n_0}}$.

{\it Case 2:} Consider now the case of a component $S$ of $G_{n_0}\cap B_{\Z}$ such that $\b S\cap\b G_{n_0}\not=\emptyset$ and denote by $h$ the regular fiber of $S$ and set $\hat{S}=\alpha(S)$.  Since $f_{\ast}(h)\not=1$ in $\pi_1N$ then the same argument as before implies that the Seifert fibration of $S$ extends to a Seifert fibration of $\hat{S}$. But since $f_{\ast}(\pi_1(S))={\Z}$ then $\hat{S}$ can be homeomorphic to a fibered  solid torus.

If  $\hat{S}$ is a solid torus $V_T$ then this means  that   $S$ has a single component $T$ in ${\rm int}(G_{n_0})$ and $\b S-T$ is adjacent  to Seifert fibered pieces in $B$. Let $S'$ be the Seifert fibered piece in $G_{n_0}$ which is adjacent to $S$ along $T$. It follows from the construction that the regular fiber of $\hat{S}'=\alpha(S')$ represented in $T$ is not free homotopic to the meridian of $\b\hat{S}=T=\b V_T$. Consider the space $\o{S}=\hat{S}'\cup_T\hat{S}$ . Thus the Seifert fibration of $\hat{S}'$ extends to a Seifert fibration of $\o{S}$. If  $\o{S}$ is not a solid torus then $\b\o{S}$ is incompressible and we have a reduction to the first case. This is true in particular when $f|S':S'\to N$ is a non-degenerate map. If $\o{S}$ is still a solid torus then we reiterate the same process. This process must stop. To see this it is sufficient to check the following

\begin{claim}\label{stop} If $G^i$ is a component of  $\hat{G}_{n_0}$ which is sent via $f_1$ into $N$ with non-zero degree then $G^i$ contains at least one non-degenerate Seifert piece.
\end{claim} 

\begin{proof}[Proof of Claim \ref{stop}]
Suppose that each Seifert piece of $G^i$ is degenerate under $f_1|G^i$. Then each Seifert piece $S$ of $G^i$ satisfies $(f_1)_{\ast}(\pi_1S)\simeq{\Z}$ by the construction. This implies, using condition (2),   that $G^i$ is a graph manifold and that the canonical tori of $G^i$ are degenerate under $f_1$. Now, since $f_1|G^i\co G^i\to N$ has nonzero degree, then using the same construction as above, one can show that there exists a Seifert fibered space $\hat{S}$, obtained from a Seifert piece $S$ in $G^i$ after Dehn filling, and a nonzero degree map $\hat{f}_1\co\hat{S}\to N$ such that $f_1|S\simeq\hat{f}_1\circ\alpha$, where $\alpha\co S\to\hat{S}$ denotes the canonical quotient map. Since $\hat{f}_1$ has nonzero degree and since $(\hat{f}_1)_{\ast}(\pi_1\hat{S})$ is cyclic this means that  $\pi_1N$ contains a cyclic finite index subgroup. This is impossible since $N$ is a non geometric closed Haken manifold. 
\end{proof}

This proves that $\hat{G}_{n_0}$ is still a Haken manifold with a torus decomposition induced from that of $M$.
 Moreover :
$${\rm Vol}(\hat{G}_{n_0})\geq{\rm deg}(f_1){\rm Vol}(N)={\rm deg}(f){\rm Vol}(N)={\rm Vol}(M)$$
and by condition (2)
$${\rm Vol}(\hat{G}_{n_0})={\rm Vol}(\hat{G}_{n_0}\cap{\cal H}(M))={\rm Vol}({\rm int}({\cal H}(M)))={\rm Vol}(M)$$
thus ${\rm Vol}(\hat{G}_{n_0})={\rm deg}(f_1){\rm Vol}(N)$.
Applying  Theorem \ref{ndeg}  to the set of non-degenerate Seifert  pieces ${\cal S}_0(\hat{G}_{n_0})$ of $\hat{G}_{n_0}$ we may assume after a homotopy supported on a regular neighborhood of  ${\cal S}_0(\hat{G}_{n_0})$ that $f_1({\cal S}_0(\hat{G}_{n_0}))\subset{\rm int}({\cal S}(N))$. Let $S$ be a degenerate Seifert piece in $\hat{G}_{n_0}$ adjacent along a canonical torus $T$ to an element $S'$ in ${\cal S}_0(\hat{G}_{n_0})$.   After a homotopy on a small regular neighborhood of $S$ we may assume, by Lemma \ref{mapping}, that $f(S)\subset{\rm int}({\cal S}(N))$. Since each component $G^i$  of $\hat{G}_{n_0}$ satisfies deg$(f_1|G^i:G^i\to N)\not=0$ (see Remark \ref{alacon}) then it contains some non-degenerate Seifert fibered pieces and thus    we may assume by repeating our argument  that $f(B_{\Z})\subset$ int$({\cal S}(N))$. Hence $f_1({\cal S}(\hat{G}_{n_0}))\subset{\rm int}({\cal S}(N))$. This ends the proof of  Claim \ref{prout} and  completes the proof of Lemma \ref{splitting}. The proof of Lemma \ref{haken} follows  directly from Lemma \ref{splitting}.
\end{proof}

\subsection{Proof of Proposition \ref{B}} 
Let $(N_i)_{i\in{\N}}$ be a sequence of non-geometric closed Haken manifolds   such that for each $i\in{\N}$ there exists a  degree-$d$ map $g_i\co M\to N_i$ with ${\rm Vol}(M)={\rm deg}(g_i){\rm Vol}(N_i)$. Throughout the proof of Proposition \ref{B} one can assume, without loss of generality, that the targets satisfy the following condition:

(III) for any $i\in{\N}$, each Seifert piece of $N_i$ has orientable orbifold base and admits an ${\bf H}^2\times{\R}$-geometry.

Indeed, let $N_i$ denote a non-geometric closed Haken manifold and let $S$ be a Seifert piece of $N_i\setminus{\mathcal T}_{N_i}$. Notice that since $N_i$ is  non-geometric then $S$ has non-empty boundary. If $S$ does not admits an ${\Hi}^2\times{\R}$-geometry this means that $\chi(S)\geq 0$, where $\chi(S)$ denotes the Euler Characteristic of the base 2-orbifold of a Seifert fibration on $S$. Since $\b S$ is non-empty and incompressible then $S$ is a Seifert fibered space over the disk with exactly two singular fibers of type $(2,1)$ which is homeomorphic to the orientable ${\S}^1$-bundle over the Moebius band.  

If $S$ is a geometric piece of $N_i$ with a Seifert fibration  over a non-orientable orbit surface then $S$ has a double cover $\t{S}$ corresponding to the orientation double cover of its orbit surface. Note that this double cover is trivial on the boundary and thus the components of $\b S$ lift to this cover. Then by taking a copy of this double cover for each component of ${\mathcal S}(N)$ that admits a Seifert fibration over a non-orientable surface, taking two copies of each component otherwise and identifying these along their torus boundary via the sewing involution between the components of $N_i\setminus{\mathcal T}_{N_i}$ (since the boundary components of each component of $N_i\setminus{\mathcal T}_{N_i}$ lift then so do the sewing involution)  we obtain a double cover $p_i\co\t{N}_i\to N_i$ satisfying condition (III) above. Note that when $S$ is a Seifert fibered piece of $N_i$  then, either

(a)  $\chi(S)<0$, $S$ has an orientable orbit space and it is covered by exactly two components $S_1, S_2$ in $\t{N}_i$ and $p_i|S_j$ is the identity, $j=1,2$, or

 (b) $\chi(S)<0$, $S$ has a non-orientable orbit space and it is covered by exactly one component $\t{S}$ in $\t{N}_i$ and $p_i|\t{S}$ is the double cover  corresponding to the orientation double cover of the orbit surface of $S$, or
 
 (c) $S$ is the orientable ${\S}^1$-bundle over the Moebius band and it is covered by a component $\t{S}$ of $\Sigma(\t{N}_i)$ homeomorphic to ${\S}^1\times{\S}^1\times I$ that can be seen as a regular neighborhood of a component of ${\mathcal T}_{\t{N}_i}$ (since $\t{N}_i$ is  non-geometric)  and $p_i|\t{S}$ is the double cover  corresponding to the orientation double cover the orbit surface of $S$. 
We have to check the following claim (notations are the same as above).

\begin{claim}\label{dom_covering}  If the family $\{N_i, i\in{\N}\}$ is infinite, up to homeomorphism, then so is  the family $\{\t{N}_i, i\in{\N}\}$.
\end{claim}  
\begin{proof} Suppose the contrary. Then we may assume, passing to a subsequence, that the family $\{\t{N}_i, i\in{\N}\}$ contains a unique element $\t{N}$ and that the family $\{N_i, i\in{\N}\}$ is infinite up to homeomorphism. First notice that the number of connected components of ${\mathcal T}_{N_i}$ is bounded by that of ${\mathcal T}_{\t{N}}$.  On the other hand, each geometric component of $N_i$ is finitely covered, via $p_i$, by a component of $\Sigma(\t{N})$. Hence, by Corollary \ref{relative}, we may assume, after passing to a subsequence, that the $N_i^{\ast}$'s are homeomorphic. For each $i\in{\N}$ we denote by $s_i$ the sewing involution that produces $N_i$ from $N_i^{\ast}$.  Let $A$ be a component of $N_i^{\ast}$ and let $Q(A)_i$ be a component of $\Sigma(\t{N})$ that covers $A$. Passing to a subsequence we may assume that $Q(A)_i$ is independant of $i\in{\N}$ and we denote it by $Q(A)$. Let $T_A$ be a component of $\b A$ and let ${\mathcal T}_A=\{U^1_A,...,U^p_A\}$ denote the components over $T_A$ in $\b Q(A)$ (again we may assume that ${\mathcal T}_A$ is independant of $i\in{\N}$). It follows from points (a), (b) and (c) of the construction that $p_i|Q(A)\co Q(A)\to A$ is either the identity or the double cover  corresponding to the orientation double cover of the orbit surface of $A$ according to whether $A$ has an orientable orbit space or not. This shows that if $l$ denotes a simple closed curve in  ${\mathcal T}_A$ then $\{(p_i)_{\ast}(l), i\in{\N}\}$ generates only one isotopy classe of curves in $T_A$. This proves that there is only one  isotopy classe of sewing involutions $s_i$ when $i\in{\N}$. Then the family  $\{N_i, i\in{\N}\}$ is finite up to homeomorphism, which gives a contradiction. 
\end{proof}  

 Since the maps   $g_i$'s are degree-$d$ maps then  by Proposition \ref{obvious} the index of $(g_i)_{\ast}^{-1}(\pi_1\t{N}_i)$ in $\pi_1M$ takes at most finitely many values. Let $\t{M}_i$ be the finite cover of $M$ corresponding to $(g_i)_{\ast}^{-1}(\pi_1\t{N}_i)$ and  let $\t{g}_i\co\t{M}_i\to\t{N}_i$ be the nonzero degree   map that covers $g_i\co M\to N_i$. Since any finitely presented group has only finitely many subgroup of given index then  there are only finitely many homeomorphisms types among $\t{M}_i$  when $i\in{\N}$.  On the other hand, if ${\rm deg}(g_i){\rm Vol}(N_i)={\rm Vol}(M)$ then ${\rm deg}(\t{g}_i){\rm Vol}(\t{N}_i)={\rm Vol}(\t{M}_i)$. Then there is no loss of generality  assuming that the targets satisfy condition (III).
 
 Now since the $N_i$'s satisfy condition (III) then the targets contains no embedded Klein bottles and one can apply Lemma \ref{haken} to the sequence of degree-$d$ maps $g_i\co M\to N_i$. Hence,  possibly after passing to a subsequence, wa can assume  that  there exists a closed Haken manifold $M_1$   which admits nonzero degree maps $f_i\co M_1\to N_i$ satisfying  properties (i), (ii) and (iii) of Lemma \ref{haken} for $i\in{\N}$.

Note that the Haken number of the $N_i$'s is bounded by that of $M_1$ and then the number of connected components of
$(N_i^{\ast},{\mathcal T}_{N_i})$ is bounded by a constant which only depends on $M_1$.  Then combining Corollary \ref{relative} with point (iii) of Lemma \ref{haken} we conclude      that there are at most finitely many topological type for $N_i^{\ast}$, when $i\in{\N}$.  This completes the proof of Proposition \ref{B}.

\section{Control of the sewing involutions of the targets}
\subsection{Statement of the Key Result for the proof of Proposition \ref{C}} 
The purpose of this section is to complete the  proof of Proposition \ref{A}. To do this it remains to prove Proposition \ref{C}. Let $d$ be a strictly positive integer.   Let $(N_i)_{i\in{\N}}$ be a sequence of weakly equivalent
non-geometric closed Haken manifolds   such that for each $i\in{\N}$ there exists a  degree-$d$ map $g_i\co M\to N_i$ with ${\rm Vol}(M)={\rm deg}(g_i){\rm Vol}(N_i)$. As in paragraph 5.4 one can assume that the $N_i$'s satisfy condition (III). Possibly after passing to a subsequence,  one can assume, using Lemma \ref{haken}, that there exists a closed Haken manifold   $M_1$ and a nonzero degree map  $f_i\co M_1\to N_i$   satisfying the following properties:

(i) ${\rm Vol}(M_1)={\rm deg}(f_i){\rm Vol}(N_i)$,

(ii) the map $f_i$ induces a finite covering between ${\mathcal H}(M_1)$ and  ${\mathcal H}(N_i)$,

(iii) for any $Q$ in $N_i^{\ast}$ each component of $(f_i)^{-1}(Q)$ is a canonical submanifold of $M_1$.

\begin{remark}\label{mini}
For convenience one require the following additional condition for point (iii) of Lemma \ref{haken}. Over all maps homotopic to $f_i\co M_1\to N_i$ satisfying point (iii) we choose always the maps such that the number of connected components of $(f_i)^{-1}(N_i^{\ast})$ is minimal. Since $N_i$ is nongeometric this implies, using Lemma \ref{mapping},  the following property: Let $B_i$ be a component of ${\mathcal S}(N_i)$ and let $W_i$ be a component of $(f_i)^{-1}(B_i)$.  Then $W_i$   contains at least one geometric piece $Q_i$ such that $(f_i)_{\ast}(t_{Q_i})\in\l h_{B_i}\r$, where $t_{Q_i}$ (resp. $h_{B_i}$) denotes the regular fiber in $Q_i$ (resp. in $B_i$). 
\end{remark}
A nonzero degree map between closed Haken manifolds satisfying points (i), (ii), (iii) and the minimality property of Remark \ref{mini} will be termed on \emph{standard form}.
 
Denote by $Q_1,...,Q_l$ the component of the $N_i^{\ast}$'s and by $T_1^k,...,T_{n_k}^k$ the boundary components of $Q_k$, for $k=1,...,l$. Each sewing involutions $s_i$, $i\in{\N}$, induces a fixed point free bijection denoted by $s_i^{\ast}$ on the set $\{T_1^k,...,T_{n_k}^k, 1\leq k\leq l\}$ by setting
$$s_i^{\ast}\co(v,i)\mapsto(w,j)\ \ {\rm if} \ \ T^v_i \ \ {\rm is\ identified\ with}\ \ T^w_j.$$
Thus, passing to a subsequence we may assume that 

(IV) for any $i, j$ in ${\N}$ then   $s_i^{\ast}=s_j^{\ast}$.	

Moreover, throughout the proof of Proposition \ref{C}, we claim that there is no loss of generality assuming that the targets satisfy the following condition:

(V) any connected component of $N_i^{\ast}$, $i\in{\N}$, has at least two boundary components.

Condition (V) comes from the following result.
\begin{lemma}\label{covering} 
Let $\{N_i\}_{i\in{\N}}$ be a sequence of weakly equivalent non-geometric Haken manifolds satisfying conditions (III) and (IV). Then there exists an integer $d>0$ such that 
for each $i\in{\N}$  there exists a $d$-fold  covering $p_i\co\t{N}_i\to N_i$ of $N_i$ such that 

(i) each component of  $\t{N}_i\setminus{\mathcal T}_{\t{N}_i}$ has at least two  boundary components,

(ii) the family $\{\t{N}_i, i\in{\N}\}$ is a sequence of weakly equivalent  Haken manifolds satisfying conditions (III) and (IV),  

(iii) if the family $\{N_i, i\in{\N}\}$ is infinite, up to homeomorphism, then so is $\{\t{N}_i, i\in{\N}\}$.
\end{lemma}
We will use the following terminology for convenience. Let ${\mathcal T}$ be a 2-manifold whose components are all tori and let $m$ be a positive integer. A covering space  $\t{\mathcal T}$ of ${\mathcal T}$ will be termed $m\times m$-\emph{characteristic} if each component of $\t{\mathcal T}$ is equivalent to the covering space of some component $T$ of ${\mathcal T}$  associated to the characteristic subgroup $H_m$ of index $m\times m$ in $\pi_1T$ (if we identify $\pi_1T$ with ${\Z}\times{\Z}$ then $H_m=m{\Z}\times m{\Z}$).
\begin{proof}[Proof of Lemma \ref{covering}]
 Since $\{N_i\}_{i\in{\N}}$ is a sequence of weakly equivalent non-geometric closed Haken manifolds then $N_i^{\ast}$ is homeomorphic to $N_j^{\ast}$ for any $i,j\in{\N}$. Then we denote by $Q_1,...,Q_l$ the component of  $N_i^{\ast}$, $i\in{\N}$. Since each $N_i$ satisfies condition (III) then using Theorems 2.4 or 3.2 of  \cite{L}, according to whether $Q_j$ is Seifert fibered or hyperbolic, we know that there is a prime $q$ such that for every $j=1,...,l$, there is a finite regular covering $p_j\co\t{Q}_j\to Q_j$ such that for any component $T$ of $\b Q_j$ then $(p_j)^{-1}(T)$ consists of more than one component and for any component $\t{T}$ of $\b\t{Q}_j$ over $T$ then  $p_j|\t{T}\co\t{T}\to T$ is the $q\times q$-characteristic covering. Denote by $\eta_j$ the degree of $p_j\co\t{Q}_j\to Q_j$. Then $(p_j)^{-1}(T)$ consists of exactly $\eta_j/q^2\geq 2$ copies of a torus. Let $m={\rm l.c.m.}(\eta_1/q^2,...,\eta_l/q^2)$. Take $t_j=m/(\eta_j/q^2)$ copies of $\t{Q}_j$, $j=1,...,l$, and glue the component of $\coprod_{j=1,l}(\coprod_{1,t_j}\t{Q}_j)$ together via lifts of the sewing involution $s_i$ of $N_i$ in the following way: let $T$ be a component of $\b Q_j$ and $T'$ be a component of $\b Q_k$ such that $T$ is identified to $T'$ in $N_i$ via $s_i|T\co T\to T'$. Note that by Condition (IV), the couple $(T,T')$ does not depend  on $i\in{\N}$. Let $\t{T}, \t{T}'$ be components of $(p_j)^{-1}(T), (p_k)^{-1}(T')$. Since $p_j|\t{T}\co\t{T}\to T$ and $p_k|\t{T}'\co\t{T}'\to T'$ are both the $q\times q$-characteristic covering then there is a sewing involution $\t{s}_i$ such that $\t{s}_i|\t{T}\co\t{T}\to\t{T}'$ covers $s_i|T\co T\to T'$. This gives a finite covering $\t{N}_i$ of $N_i$ satisfying properties (i) and (ii). To check property (iii) it is sufficient to apply the same kind of arguments as in the proof of Claim \ref{dom_covering}.
\end{proof}

 Then, to complete the proof of Proposition \ref{C} we first   state  the following technical key result which shows that the sewing involution of the domain "fix", in a certain sense, the sewing involution of the targets. This result combined with Lemma \ref{final} (section 6.3) ensures the finiteness of the equivalence classes of the sewing involutions of the targets. 
\begin{lemma}[Gluing Lemma]\label{end} Let $M_1$ be a closed Haken manifold and let $\{N_i, i\in{\N}\}$ be a sequence of weakly equivalent non-geometric closed Haken manifolds satisfying conditions (III), (IV), and (V) such that there exist nonzero degree maps $f_i\co M_1\to N_i$ on standard form satisfying  ${\rm Vol}(M_1)={\rm deg}(f_i){\rm Vol}(N_i)$. Let $A$ and $B$ be two components of $(N_i^{\ast},s_i)$ such that
$s_i$ connects a component $T_A$ of $\b A$ with a component $T_B$ of $\b B$ for any $i\in{\N}$. Denote by $T$ the component of ${\mathcal T}(N_i)$ obtained by identifying $T_A$ with $T_B$ via $s_i$. Then, possibly after passing to a subsequence, the following properties hold:

(i) if $A$ and $B$ are both hyperbolic pieces of $N_i$ then the maps $\{s_i|T_A\co{T_A}\to{T_B}\}_{i\in{\N}}$ are in the same isotopy class,

(ii) if $A$ and $B$ are both Seifert pieces of $N_i$ then there exists two elements $(a,b)$ and $(c,d)$ of ${\mathcal P}_{\ast}$, which depend on $T$, such that $bd\not=0$ and a sequence $\{\delta_A^i\}_{i\in{\N}}$ (resp. $\{\delta_B^i\}_{i\in{\N}}$) of sections for $T_A$ (resp. for $T_B$), with respect to the Seifert fibration of $A$ (resp. of $B$)  such that   
$$(s_i)_{\ast}(h_{A})=h^{a}_{B}(\delta_B^i)^{b}\ \ \ {\rm and}\ \ \    (s_i)_{\ast}(h_{B})=h^{c}_{A}(\delta_A^{i})^{d}$$
where $h_A$ (resp. $h_B$) denotes the regular fiber of $A$ (of $B$ resp.)

(iii) if $B$ is a Seifert piece and if $A$ is hyperbolic then there exists a basis $(\lambda_A,\mu_A)$ of $\pi_1T_A$ and a sequence $\{\delta_i,i\in{\N}\}$ of sections of $T_B$ with respect to the Seifert fibration of $B$ such such that     
$$(s_i)_{\ast}(\lambda_A)=h_{B}^{\pm 1}\ \ \ {\rm and}\ \ \  (s_i)_{\ast}(\mu_A)=\delta_{i}$$ for any $i\in{\N}$.

\end{lemma}

\subsection{Proof of Lemma \ref{end}}
Denote by $T$ the component of ${\mathcal T}_{N_i}$ obtained by sewing $A$ and $B$ via $s_i|T_A\co T_A\to T_B$.

\emph{Case 1: $A$ and $B$ are hyperbolic manifolds}. There exists two components $H_A$ and $H_B$ of ${\mathcal H}(M_1)$ (independant of $i\in{\N}$, possibly after passing to a subsequence) adjacent along some components of $f_i^{-1}(T)$ such that $f_i|H_A\co H_A\to A$ and $f_i|H_B\co H_B\to B$ are nonzero degree maps (actually these maps are finite coverings). Possibly after passing to a subsequence, we may assume that $(f_i|H_A)^{-1}(T_A)$ is independant of $i\in{\N}$. Denote by $\{U^1_A,...,U^p_A\}$ the component of $(f_i|H_A)^{-1}(T_A)$.  Then to prove  Lemma \ref{end} in the first case it is sufficient to show the following result.
\begin{claim}\label{endhh}
 Let $l$ be a  simple closed curve  in $U^1_A$ and let $l_T^i$ be the simple closed curve in $T_A$ such that   $(f_i)_{\ast}([l])\in\l[l_T^i]\r$. Then the set of curves $\{l_T^i, i\in{\N}\}$ generates at most finitely many isotopy classes of curves in $T_A$.
\end{claim}
\begin{figure}[h]
\centerline{\input{chirurgie.pstex_t}}
\end{figure}
\begin{proof}[Proof of Claim \ref{endhh}]
  Consider the 3-manifold $A_i$ obtained after performing a Dehn filling on $A$ by identifying the meridian of a solid torus with $l_T^i$ and denote by $r_i\co A\to A_i$ the canonical quotient map.   For each $v\in\{2,...,p\}$ there exists a  simple closed curve $l^i_v\in U^v_A$ such that $(r_i\circ f_i)(l^i_v)$ is nul homotopic in $\pi_1 A_i$. Denote by $H_A^l$ the 3-manifold obtained from $H_A$ by performing a Dehn filling along $U^1_A$ identifying the meridian of a solid torus with $l$ and denote by $H_A^i$ the 3-manifold obtained from $H_A^l$ by gluing $p-1$ solid tori along $U_A^2\cup...\cup U_A^p$ by identifying each meridian with  $l^i_v$ when $v\in\{2,...,p\}$. Note that when $p=1$ then $H_A^l=H_A^i$ is independant of $i\in{\N}$. Using the Mayer Vietoris exact sequence it is easily checked that $f_i|H_A\co H_A\to A$ induces a proper nonzero degree map   $\hat{f}_i\co H_A^i\to A_i$  such that ${\rm deg}(\hat{f}_i)={\rm deg}(f_i|H_A)$. Assume that the curves $l_T^i$ generate   infinitely many isotopy classes of curves in $T_A$. Thus when the lenght of the curve $l_T^i$ is sufficiently large then the formulae established in \cite{NZ} implies that 
  $${\rm Vol}(A_i)\approx{\rm Vol}(A)-\pi^2\frac{{\mathcal A}(T_A)}{{\rm lenght}(l^i_T)}$$ 
  where ${\mathcal A}(T_A)$ denotes the area of the torus $T_A$ with respect to the Euclidean structure induced by the complete hyperbolic structure on ${\rm int}(A)$ and where ${\rm lenght}(l^i_T)$ is the lenght of the curve $l^i_T$ on the torus $T_A$ with respect to this Euclidean structure.   
   Then we may assume, passing to a subsequence, that $\{{\rm Vol}(A_i), i\in{\N}\}$ is a striclty increasing sequence such that
 $$\lim_{i\to\infty}{\rm Vol}(A_i)={\rm Vol}(A)\ \ \ \ \ (\star)$$  
 and that the $A_i$'s 	are complete finite volume hyperbolic 3-manifolds by the Hyperbolic Surgery Theorem of W. P.  Thurston, \cite{Th6}.
  
 Assume first that $p=1$ (i.e. $(f_i|H_A)^{-1}(T_A)=U^1_A$ is connected). Then the latter equality implies that $H_A^l$ dominates infinitely many hyperbolic manifolds (these manifolds can be distinguished by their volume) which contradicts Corollary \ref{relative}. When $p\geq 2$  then  consider the following sequence of inequalities which hold for any $i\in{\N}$:
$${\rm Vol}(H_A)>{\rm Vol}(H_A^l)>{\rm Vol}(H_A^i)\geq{\rm deg}(\hat{f}_i){\rm Vol}(A_i)={\rm deg}(f_i|H_A){\rm Vol}(A_i)$$ 
Since  $f_i|H_A\co H_A\to A$ is a finite covering between hyperbolic manifolds then ${\rm deg}(f_i|H_A)$ is a constant equal to ${\rm Vol}(H_A)/{\rm Vol}(A)$. Then using equality $(\star)$ we have $$\lim_{i\to\infty}{\rm deg}(f_i|H_A){\rm Vol}(A_i)={\rm deg}(f_i|H_A){\rm Vol}(A)={\rm Vol}(H_A)$$ Hence we get a contradiction. This proves the claim.
  
\end{proof}  
\begin{proof}[End of proof of Lemma \ref{end} Point (i)]
Choose a basis $(\t{\lambda}_A,\t{\mu}_A)$ of $\pi_1U^1_A$. Since $f_i|H_A\co H_A\to A$ is a covering  then $\t{\lambda}_A$ (resp. $\t{\mu}_A$) can be identified to an element $\lambda_A$ (resp. $\mu_A$) of $\pi_1T_A$.  Let $c\in\pi_1T_A$ be a primitive element. Since $f_i|H_A\co H_A\to A$ is a finite covering then there exists an integer $n$ such that $c^n$ is a primitive element in $\pi_1U^1_A$. The sewing involution $\sigma_1$ in $M_1$ identifies $c^n$ with a primitive element $\t{l}$ in $U^1_B$, where $U^1_B$ is the component of $\b H_B$ such that $U^1_B=\sigma_1(U^1_A)$. Using Claim \ref{endhh} then, after passing to a subsequence, we may assume that there exists a simple closed curve $l_T$ in $T_B$ such that $(f_i|B)_{\ast}([\t{l}])\in\l[l_T]\r$, when $i\in{\N}$. This proves, possibly  after passing to a subsequence, that $(s_i)_{\ast}(c)=[l_T]^{\pm 1}$ for any $i\in{\N}$, where $c$ is an arbitrary element of $\pi_1T_A$. This completes the proof of Lemma \ref{end} point (i).
\end{proof}
\begin{figure}[h]
\centerline{\input{sewing.pstex_t}}
\end{figure}
\emph{Case 2: $A$ and $B$ are Seifert fibered spaces.} 
 Fix a component $Q_A^i$  of the preimage of $A$ such that $f_i|Q_A^i\co Q_A^i\to A$ has nonzero degree and let $Q_B^i$ denote the components of the preimage of $B$ which are adjacent to $Q_A^i$ along $(f_i|Q_A^i)^{-1}(T_A)$. Passing to a subsequence we may assume that $Q_A^i$ and $Q_B^i$ are independant of $i$ and we denote them by $Q_A$ and $Q_B$. Note that, by Remark \ref{debile}, $Q_A$ and $Q_B$ are graph manifolds. 

\subsubsection{First Step}
For convenience, we perform the following modification on $Q_A$. Let $S^j_B$ be a Seifert piece of $Q_B$ adjacent to $Q_A$ along a component of $(f_i|Q_A)^{-1}(T_A)$ such that $(f_i)_{\ast}(h^j_B)\not\in\l h_B\r$, where $h^j_B$ denotes the regular fiber of $S^j_B$ and $h_B$ is the regular fiber of $B$.  Then by Lemma \ref{mapping} one can perturb slightly $f_i$ by a homotopy, which is constant outside of a regular neighborhood of $S^j_B$,  so that $Q_A^{i,j}=Q_A\cup S^j_B$ is a component of $(f_i)^{-1}(A)$. We do that for any  Seifert  piece $S^j_B$ of $Q_B$ adjacent to $Q_A$ along a component of $(f_i|Q_A)^{-1}(T_A)$ such that $(f_i)_{\ast}(h^j_B)\not\in\l h_B\r$.  Denote  by $Q_A^{i,{\rm new}}$ (resp. $Q^{i,{\rm new}}_B$) the \emph{new preimage} of $A$ (resp. of $B$). After repeating this process a finite number of times we may assume that each Seifert piece  $W$ of $Q^{i,{\rm new}}_B$ adjacent to a Seifert piece of $Q_A^{i,{\rm new}}$ along a component of $(f_i|Q_A^{i,{\rm new}})^{-1}(T_A)$ satisfies: $f_{\ast}(h_W)\in\l h_B\r$, where $h_W$ denotes the homotopy class of the regular fiber in $W$. Passing to a subsequence we may assume that $Q_A^{i,{\rm new}}$  and $Q^{i,{\rm new}}_B$ are independant of $i\in{\N}$ and we denote them by $Q_A^{{\rm new}}$ and $Q^{{\rm new}}_B$.  Note that $Q^{{\rm new}}_B\not=\emptyset$ by Remark \ref{mini}. Denote by ${{\mathcal T}}_A^{\rm new}=\{U_A^1,...,U_A^l\}$ the components of $(f_i|Q_A^{{\rm new}})^{-1}(T_A)$ and by ${{\mathcal T}}_B^{\rm new}=\{U_B^1,...,U_B^l\}$ the components of  $\b Q^{{\rm new}}_B$ adjacent to  ${{\mathcal T}}_A^{\rm new}$. For each $i\in\{1,...,l\}$, choose a simple closed curve $t_B^i$ in $U_B^i$ which represents the regular fiber of the Seifert fibered space in ${{\mathcal T}}_B^{\rm new}$ containing $U_B^i$. By construction we know that $(f_i)_{\ast}([t_B^i])\in\l [h_B]\r$. Using the sewing involution $\sigma_1$ of $M_1$   the family of curves $\{t_B^1,...,t_B^l\}$ define a family of curves $\{c_A^1,...,c_A^l\}$ in ${{\mathcal T}}_A^{\rm new}$ defined by $c^j_A=\sigma_1(t^j_B)$. It follows from our construction that for any $i\in{\N}$ there exists a simple closed curve $l_T^i$ in $T_A$ such that   $(f_i)_{\ast}([c_A^v])\in\l[l_T^i]\r$ for $v=1,...,l$ and $i\in{\N}$.

\subsubsection{Second Step}
 Consider the 3-manifold $A_i$ obtained after performing a Dehn filling on $A$ by identifying the meridian of a solid torus with $l_T^i$ and denote by $r_i\co A\to A_i$ the canonical quotient map.    Denote by ${\mathcal D}(Q_A^{{\rm new}})$ the 3-manifold obtained from $Q_A^{{\rm new}}$ after gluing $l$ solid tori along ${{\mathcal T}}_A^{\rm new}$ by identifying each meridian with  $c_A^v$ when $v\in\{1,...,l\}$. This gives proper  nonzero degree maps $\hat{f}_i\co {\mathcal D}(Q_A^{{\rm new}})\to A_i$ with $\b A_i\not=\emptyset$ and $\b{\mathcal D}(Q_A^{{\rm new}})\not=\emptyset$ by condition (V).

Fix a section $\delta_0$ on $T_A$ with respect to the Seifert fibration of $A$.  Let $(a_i,b_i)$ be a sequence of coprime integers such that $$l^i_T=h_A^{a_i}\delta^{b_i}_{0}, \ \ {\rm for}\ i\in{\N}$$
We first claim that the sequence $\{b_i,i\in{\N}\}$ is finite. Indeed, if
 the sequence $\{b_i,i\in{\N}\}$ is infinite then  we get infinitely many pairwise non-homeomorphic Seifert fibered spaces $\{A_i,i\in{\N}\}$ properly dominated by ${\mathcal D}(Q_A^{{\rm new}})$. These Seifert fibered spaces can be distinguished, for example, by the order of the exceptional fiber generated by performing the Dehn fillings along $T_A$.    This gives a contradiction with Corollary \ref{relative} since the $A_i$ have non-empty boundary with geometry ${\bf H}^2\times{\R}$. Thus from now one one can assume, passing to a subsequence, that $b_i$ is a constant denoted by $b$. 
 
 Consider now the sequence $\{a_i,i\in{\N}\}$. We know that $A_i$ is homeomorphic to $A_j$ if and only if $a_i=a_j\ {\rm mod}(b)$. Then using the same argument as above,  we may assume, possibly after passing to a subsequence, that there exists an integer $a$ and a sequence of integers of  $\{m_i, i\in{\N}\}$ such that $a_i=a+m_ib$, for any $i\in{\N}$. Hence we get $$l^i_T=h_A^{a+m_ib}\delta_0^b=h_A^a\delta_i^b, \ \ {\rm for}\ i\in{\N}$$ where   $\delta_i$ is the section of $T_A$ defined by $c_i=\delta_{0}h_A^{m_i}$. Thus $(s_i)_{\ast}(h_B)=h_A^a\delta_i^b$, $i\in{\N}$. Notice that $b\not=0$ by the minimality property of the geometric decomposition for Haken manifolds.  The proof of point (ii) of Lemma \ref{end} follows by permuting the role of $Q_A$ and $Q_B$. Indeed it is sufficient to choose a component of $(f_i)^{-1}(B)$ which dominates $B$ and to proceed in the same way as above.

\emph{Case 3: $B$ is Seifert and $A$ is a hyperbolic manifold.}
Let $Q_B$ be a component of the preimage of $B$ such that $f_i|Q_B\co Q_B\to B$ has nonzero degree and let $Q_A$ be the preimage (not necessarily connected) of $A$ which is adjacent to $Q_B$ along ${\mathcal T}_B=(f_i|Q_B)^{-1}(T_B)$. Denote by ${\mathcal T}_A$ the components of $\b Q_A$ identified with ${\mathcal T}_B$ in $M_1$ (as in the paragraphs above we can assume that $Q_A, Q_B, {\mathcal T}_A$ and  ${\mathcal T}_B$ are independant of $i\in{\N}$). Note that it follows from Remark \ref{debile} that $Q_A$ is a disjoint union of hyperbolic manifolds. Moreover we may assume, passing to a subsequence that $f_i|Q_A\co Q_A\to A$ are equivalent coverings, when $i\in{\N}$. We first state the following technical result.
\begin{figure}[h]
\centerline{\input{mixte.pstex_t}}
\end{figure}
\begin{claim}\label{endsh}
Let $c$ be a simple closed curve in ${\mathcal T}_B$  and let $h_B^i$  be the simple closed curve in $T_B$  such that $(f_i|Q_B)_{\ast}([c])\in\l[h^i_B]\r$. Then, possibly after passing to a subsequence, the following properties hold:
there exists an element $(a,b)$ in ${\mathcal P}_{\ast}$ and a sequence of sections $\{\delta_i,i\in{\N}\}$ of $T_B$ with respect to the Seifert fibration on $B$ such that 
$$h_B^i=h_B^a\delta_i^b, \ i\in{\N}$$ where $h_B$ denotes the homotopy class of the regular fiber on $B$.
\end{claim}
\begin{proof}[Proof of Claim \ref{endsh}]
Denote by $\{T_A^1,...,T_A^l\}$ (resp. $\{T_B^1,...,T_B^l\}$) the components of 
${\mathcal T}_A$ (resp. ${\mathcal T}_B$)  in such a way that $T^v_A$ is identified with $T^v_B$ in $M_1$ when $v=1,...,l$.
Let $c=c_B^1$ be a simple closed curve on $T_B^1$ and let $h^1_B$ denotes the simple closed curve in $T_B$ such that $(f_1|Q_B)_{\ast}(c_B^1)\in\l[h^1_B]\r$. Complete this curve by a system of curves  $\{c_B^1,...,c_B^l\}$ where $c_B^v$ is a simple closed curve in $T_B^v$ such that   $(f_1|Q_B)_{\ast}(c_B^v)\in\l[h^1_B]\r$ for $v=1,...,l$. 

Note that since $Q_B$ is adjacent along ${\mathcal T}_B$ to hyperbolic pieces then $f_i|{\mathcal T}_B$ is $\pi_1$-injective and thus  $(f_1|Q_B)_{\ast}(c_B^v)$ is non-trivial (for $v=1,...,l$).  The sewing involution $\sigma_1$ of $M_1$ allows to asociate to the family $\{c_B^1,...,c_B^l\}$ a system of curves $\{c_A^1,...,c_A^l\}$ in ${\mathcal T}_A$ such that $c_A^v$ is a simple closed curve in $T_A^v$ and $c_A^v=\sigma_1(c_B^v)$ for $v=1,...,l$. On the other hand, denote by $h^1_A$ the simple closed curve in $T_A$ defined by $h_A^1=s_1(h_B^1)$. Thus it follows from our construction that $(f_1|Q_A)_{\ast}(c_A^v)$ are non-trivial elements of the cyclic group generated by $[h_A^1]$, for $v=1,...,l$. Since the sequence of maps  $\{(f_i|Q_A)\co Q_A\to A,i\in{\N}\}$ are equivalents finite coverings then the family $\{(f_i|Q_A)_{\ast}(c_A^1),...,(f_i|Q_A)_{\ast}(c_A^l)\}$ is always contained in an infinite cyclic subgroup of $\pi_1T_A$ when $i\in{\N}$. Then, by construction, the family  $\{(f_i|Q_B)_{\ast}(c_B^1),...,(f_i|Q_B)_{\ast}(c_B^l)\}$ satisfies the same property and thus for any $i\in{\N}$ there exists  a simple closed curve $h_B^i$ in $T_B$ such that $(f_i|Q_B)_{\ast}(c_B^v)\in\l[h_B^i]\r$ for $v=1,...,l$ and $i\in{\N}$. To complete the proof of Claim \ref{endsh} it is sufficient to apply the same arguments as in paragraph 5.2.2.  

\end{proof}
\begin{proof}[End of Proof of Lemma \ref{end} point (iii)]
Let $(\t{\lambda}_A,\t{\mu}_A)$ be a basis of $\pi_1T^1_A$. One can define a basis  $(\t{\lambda}_B,\t{\mu}_B)$ of $\pi_1T^1_B$ with $\t{\lambda}_B=\sigma_1(\t{\lambda}_A)$ and $\t{\mu}_B=\sigma_1(\t{\mu}_A)$ (where $\sigma_1$ is the sewing involution of $M_1$). By Claim \ref{endhh}, possibly after passing to a subsequence, one can assume that there exist  two primitive curves  $h_A$ and $k_A$ in $T_A$ such that $$(f_i|Q_A)_{\ast}(\t{\lambda}_A)\in\l[h_A]\r\ {\rm and}\  (f_i|Q_A)_{\ast}(\t{\mu}_A)\in\l[k_A]\r\ {\rm for\ any}\ i\in{\N}.$$ 
Note that since $f_i|Q_A\co Q_A\to A$ is a finite covering then $[h_A],[k_A]$ are algebraically free in $\pi_1T_A$ and the group generated by $(h_A,k_A)$ is a  finite index subgroup  of $\pi_1T_A$. On the other hand, using Claim \ref{endsh}, we know that there exist two elements $(a,b)$ and $(c,d)$ of ${\mathcal P}_{\ast}$  and two sequences of sections $\{\delta_i,i\in{\N}\}$, $\{d_i,i\in{\N}\}$ of $T_B$ (with respect to the Seifert fibration of $B$) such that  
$$(f_i|Q_B)_{\ast}(\t{\lambda}_B)\in\l h_B^a\delta_i^b\r\ {\rm and}\   (f_i|Q_B)_{\ast}(\t{\mu}_B)\in\l h_B^cd_i^d\r$$
 Let $\{n_i,i\in{\N}\}$ denote the sequence of integers such that $d_i=\delta_ih_B^{n_i}$ for $i\in{\N}$. With this notation we get  $(f_i|Q_B)_{\ast}(\t{\mu}_B)\in\l h_B^{c+dn_i}\delta_i^d\r$. Note that  $f_i|{\mathcal T}_B\co{\mathcal T}_B\to T_B$ is a finite covering of constant degree when $i\in{\N}$. Indeed since $f_i|Q_A\co Q_A\to A$ is a finite covering between hyperbolic 3-manifolds then ${\rm deg}(f_i|Q_A)$ is a constant equal to ${\rm Vol}(Q_A)/{\rm Vol}(A)$ and thus the degree of the maps $f_i|{\mathcal T}_A\co {\mathcal T}_A\to T_A$ is constant. Thus so is the degree of the maps  $f_i|{\mathcal T}_B\co{\mathcal T}_B\to T_B$. 
 
 Then, possibly after passing to a subsequence, we may assume that $\Delta_i=ad-bc-bdn_i$ does not vanish and is constant when $i\in{\N}$.   It follows from the construction that $$(s_i)_{\ast}(h_A)=h_B^a\delta_i^b \ \ {\rm and} \ \ (s_i)_{\ast}(k_A)=h_B^{c+dn_i}\delta_i^d$$

\emph{Case 1}: Assume first that $bd=0$. Then either $b=0$ and $(s_i)_{\ast}(h_A)=h_B^a$ where $a=\pm 1$ or $d=0$ and $(s_i)_{\ast}(k_A)=h_B^{c+dn_i}$ where $c+dn_i=c=\pm 1$. Then we set $\lambda_A=h_A$ or $\lambda_A=k_A$ depending on whether $b=0$ or $d=0$. Choose a primitive curve $\mu_A$ in $T_A$ in such a way that $(\lambda_A,\mu_A)$ is basis of $\pi_1T_A$. Since the subgroup of $\pi_1T_A$ generated by $h_A,k_A$ has finite index in $\pi_1T_A$  then there exists $n\in{\N}$ such that $\mu_A^n\in\l h_A,k_A\r$. Let $x,y$ be two integers such that $\mu_A^n=h_A^x k_A^y$.  Thus we get  $$(s_i)_{\ast}(\lambda_A)=h_B^{\pm 1}\ {\rm and}\ (s_i)_{\ast}(\mu_A^n)=h_B^{ax+yc+ydn_i}\delta_i^{bx+dy}$$

\emph{Subcase 1.1}: If $b=0$ then $d\not=0$, otherwise $\Delta_i=0$ which is impossible, and $(s_i)_{\ast}(\mu_A^n)=h_B^{ax+yc}d_i^{dy}$ where $d_i$ is a section of $T_B$ defined by $d_i=\delta_ih_B^{n_i}$. Notice that in this case then necessarily $y\not=0$. Thus $(s_i)_{\ast}(\mu_A)=h_B^{v}d_i^{w}$ where $(v,w)$ are coprime. Since $(s_i)_{\ast}(\lambda_A)=h_B^{\pm 1}$, and since $(\lambda_A,\mu_A)$ is a basis of $\pi_1T_A$ then $w=\pm 1$.

\emph{Subcase 1.2} : If $d=0$ then $b\not=0$ and thus $(s_i)_{\ast}(\mu_A^n)=h_B^{ax+yc}\delta_i^{bx}$ with $x\not=0$. Then there exists a couple $(v',w')$ of coprime integers such that $(s_i)_{\ast}(\mu_A)=h_B^{v'}\delta_i^{w'}$ where $w'=\pm 1$.

\emph{Case 2} : Assume now that $bd\not=0$. Note that since $\Delta_i$ is constant when $i\in{\N}$ then in this case the sequence $\{n_i,i\in{\N}\}$ is necessarily constant. Hence we denote $n_i$ by $n_0$ (when $i\in{\N}$). Consider the element of $\pi_1T_A$ given by $\o{\lambda}_A=h_A^dk_A^{-b}$. Then $(s_i)_{\ast}(\o{\lambda}_A)=h_B^{\Delta_i}$. Hence if $\lambda_A$ denotes the primitive element of $\pi_1T_A$ such that $\o{\lambda}_A\in\l \lambda_A\r$ then $(s_i)_{\ast}(\lambda_A)=h_B^{\pm 1}$. Choose a primitive element $\mu_A$ in $\pi_1T_A$ in such a way that $(\lambda_A,\mu_A)$ is a basis of $\pi_1T_A$. As in \emph{Case 1} we know that there exist integers $n,x$ and $y$ such that  $\mu_A^n=h_A^x k_A^y$. Then $$(s_i)_{\ast}(\mu_A^n)=h_B^{ax+yc+ydn_0}\delta_i^{bx+dy}$$  Since $(\lambda_A,\mu_A)$ is a basis of $\pi_1T_A$, since $s_i$ induces an isomorphism between $\pi_1T_A$ and $\pi_1T_B$ and since $(s_i)_{\ast}(\lambda_A)=h_B^{\pm 1}$ then there exists an integer $v$ such that  $(s_i)_{\ast}(\mu_A)=h_B^{v}\delta_i^{\pm 1}$.  
\end{proof}
 This ends the proof of Lemma \ref{end}.
 \subsection{End of proof of Proposition \ref{C}}
 To complete the proof of Proposition \ref{C} it remains to check that if $\{N_i,i\in{\N}\}$ denotes a sequence of weakly equivalent closed Haken manifolds satisfying the hypothesis of Lemma \ref{end} then $N_i$ is homeomorphic to $N_j$ for any $i,j$. More precisely,  the main purpose of this section is to state the following result.
 \begin{lemma}\label{final}
 Let $N_1$ and $N_2$ be two weakly equivalent non-geometric closed Haken manifolds  satisfying conditions (III), (IV) and (V). Denote by $s_1$, $s_2$ the sewing involution of $N_1$ and $N_2$. Suppose that for any pair of   components $A$ and $B$  of $N_i^{\ast}$ such that $s_i$ connects a component $T_A$ of $\b A$ with a component $T_B$ of $\b B$ then $s_i|T_A\co T_A\to T_B$ satisfies the conditions (i), (ii) or (iii) of Lemma \ref{end} depending on whether $A$ (resp. $B$) is a hyperbolic or a Seifert fibered manifold, $i=1,2$. This means that 
 
 (i) if $A$ and $B$ are hyperbolic manifolds then $s_1|T_A\co T_A\to T_B$ and $s_2|T_A\co T_A\to T_B$ are isotopic,
 
 (ii) if $A$ and $B$ are Seifert fibered spaces then there exist two elements $(a,b)$ and $(c,d)$ of ${\mathcal P}_{\ast}$ (which depend on the canonical torus $T$ obtained by gluing $T_A$ with $T_B$)  with $bd\not=0$  and sections $\delta_A^1, \delta_A^2$ for $T_A$ and $\delta_B^1, \delta_B^2$ for $T_B$, with respect to the Seifert fibrations of $A$ and $B$ respectively such that: 
 $$(s_1)_{\ast}(h_A)=h_B^a(\delta_B^1)^b \ \ {\rm and} \ \   (s_1)_{\ast}(h_B)=h_A^c(\delta_A^1)^d$$
$$(s_2)_{\ast}(h_A)=h_B^a(\delta_B^2)^b \ \ {\rm and} \ \   (s_2)_{\ast}(h_B)=h_A^c(\delta_A^2)^d$$ 
where $h_A$ and $h_B$ denotes the regular fiber of $A$ and $B$,

(iii) if $A$ is hyperbolic and $B$ is Seifert then  there exists a basis $(\lambda_A,\mu_A)$ of $\pi_1T_A$, an integer $a\in{\Z}$, which depends on $T$, and a section $\delta_B$ of $T_B$ with respect to the Seifert fibration of $B$ such that 
$(s_1)_{\ast}(\lambda_A)=h_B$, $(s_1)_{\ast}(\mu_A)=\delta_B$,  $(s_2)_{\ast}(\lambda_A)=h_B$ and $(s_2)_{\ast}(\mu_A)=h_B^a\delta_B$.

 Then  $N_1$ is homeomorphic to $N_2$.

 \end{lemma}
\begin{proof}
Let $A$ and $B$ be two canonical submanifolds of $N_i$ ($i=1,2$) such that $s_i$ connects a component $T_A$ of $\b A$ with a  component $T_B$ of $\b B$. Denote by $W_i$ the submanifold of $N_i$ obtained by gluing $A$ and $B$ identifying $T_A$ with $T_B$ via $s_i$. Then to prove the lemma it is sufficient to show that  $W_1$ is homeomorphic to $W_2$, for any canonical submanifolds $A$ and $B$ of $N_i$.  

Let $V_A$ and $V_B$ denote the geometric pieces of $A$ and $B$ adjacent to $T_A$ and $T_B$ respectively.
 
 \emph{Case 1: $V_A$ and $V_B$ are hyperbolic manifolds.} In this case, it follows from hypothesis (i) of the lemma that $s_1|T_A\co T_A\to T_B$ and $s_2|T_A\co T_A\to T_B$ are in the same isotopy class. This means that if $h$ denotes the identity map in $A\coprod B$ then $h\circ s_1$ and $s_2\circ h$ are isotopic.  This proves that $W_1\simeq W_2$.
 
 \emph{Case 2: $V_A$ and $V_B$ are Seifert fibered  manifolds.}
Then by point (ii) of the lemma, we know that there exist two sections $\delta^1_A, \delta^2_A$ of $T_A$ and $\delta^1_B, \delta^2_B$ of $T_B$ with respect to the Seifert fibration of $V_A$ and $V_B$ respectively such that
 $$(\ast)\ \ \ \ \ (s_i)_{\ast}(h_A)=h_B^a(\delta_B^i)^b\ \ {\rm and} \ \ (s_i)_{\ast}(h_B)=h_A^c(\delta_A^i)^d\ \ {\rm with} \ \ bd\not=0$$
 for $i=1,2$ and where $h_A$ and $h_B$ denote the homotopy class of the regular fiber of $V_A$ and $V_B$ respectively. 
 Notice that, since $bd\not=0$ then the equations $(\ast)$ (for $i=1,2$) determin the sewing involutions  $s_i|T_A\co T_A\to T_B$ and thus the manifolds $W_1$ and $W_2$.
Indeed there exists a unique element $(\alpha,\beta)$  of ${\mathcal P}_{\ast}$ such that 
$$(s_1)_{\ast}(\delta^1_B)=h_A^{\alpha}(\delta^1_A)^{\beta} \ {\rm and} \ \   
 (s_2)_{\ast}(\delta^2_B)=h_A^{\alpha}(\delta^2_A)^{\beta}$$
 This element is given by the equations $ac+\alpha b=1$ and $ad+\beta b=0$. 
Denote by $V_A'$, $V_B'$, resp. $A'$ and $B'$, the space obtained from $V_A$, $V_B$, resp. $A$ and $B$, after removing a regular neighborhood $V(h_A)$ and $V(h_B)$ of a regular fiber in $V_A$ and $V_B$. We write $\b V_A'=T_A\cup T_A'\cup T_A^1\cup...\cup T_A^p$ and  $\b V_B'=T_B\cup T_B'\cup T_B^1\cup...\cup T_B^q$ where $T_A'=\b V(h_A)$,  $T_B'=\b V(h_B)$ and $T_A^1\cup...\cup T_A^p$ denote the components of $\b V_A\setminus T_A$ and $T_B^1\cup...\cup T_B^q$ denote the components of $\b V_B\setminus T_B$.   For each component $T_A^i$, $T_B^j$ $i=1,...,p$, $j=1,...,q$, we fix a section $d_A^i$, $d_B^j$ of $T_A^i$ and $T_B^j$ respectively with respect to the Seifert fibration of $V_A'$ and $V_B'$ in such a way that $\pi_1V_A'$ and $\pi_1V_B'$ admit the following presentation: 
$$\pi_1V_A'=\l a_1,b_1,...,a_{g_A},b_{g_A},q_1,...,q_{r_A},\delta^1_A,\delta_A',d_A^1,...,d_A^p,h_A:$$ $$a_ih_Aa_i^{-1}=h_A, b_ih_Ab_i^{-1}=h_A, q_ih_Aq_i^{-1}=h_A, \delta^1_A h_A(\delta^1_A)^{-1}=h_A, \delta'_Ah_A(\delta'_A)^{-1}=h_A,$$ 
 $$d_A^ih_A(d_A^i)^{-1}=h_A, q_j^{\alpha_j}=h_A^{\beta_j}, [a_1,b_1]...[a_{g_A},b_{g_A}]q_1...q_{r_A}\delta^1_A\delta_A'.d_A^1...d_A^p=1\r$$
where $\delta'_A$ is the meridian of $V(h_A)$ and
 $$\pi_1V_B'=\l a_1,b_1,...,a_{g_B},b_{g_B},q_1,...,q_{r_B},\delta^1_B,\delta_B',d_B^1,...,d_B^q,h_B:$$ $$a_ih_Ba_i^{-1}=h_B, b_ih_Bb_i^{-1}=h_B, q_ih_Bq_i^{-1}=h_B, \delta^1_B h_B(\delta^1_B)^{-1}=h_B, \delta'_Bh_B(\delta'_B)^{-1}=h_B,$$ 
 $$d_B^ih_B(d_B^i)^{-1}=h_B, q_j^{\gamma_j}=h_B^{\lambda_j}, [a_1,b_1]...[a_{g_B},b_{g_B}]q_1...q_{r_B}\delta^1_B\delta_B'.d_B^1...d_B^q=1\r$$
 where $\delta'_B$ is the meridian of $V(h_B)$.
 
 We construct a homeomorphism $\eta\co V_A'\coprod V_B'\to V_A'\coprod V_B'$ by setting:
$$\eta(\delta^1_A)=\delta^2_A, \eta(d^i_A)=d^i_A, \eta(a_i)=a_i, \eta(b_i)=b_i, \eta(q_i)=q_i\ \ {\rm and}\ \  \eta(\delta'_A)=\delta'_Ah_A^u$$
where $u$ denotes the integer defined by $h_A^u\delta^2_A=\delta^1_A$ and 
$$\eta(\delta^1_B)=\delta^2_B, \eta(d^i_B)=d^i_B, \eta(a_i)=a_i, \eta(b_i)=b_i, \eta(q_i)=q_i\ \ {\rm and}\ \  \eta(\delta'_B)=\delta'_Bh_B^v$$
where $v$ denotes the integer defined by $h_B^v\delta^2_B=\delta^1_B$. 

Extend this homeomorphism by the identity on $(A\setminus V_A)\coprod(B\setminus V_B)$. We denote by $W'_i$ the space obtained after gluing $A'$ with $B'$ by identifying $T_A$ and $T_B$ via $s_i$, $i=1,2$. Then it follows directly from our construction that $W'_1$   is homeomorphic to  $W'_2$. Note that the homeomorphism $$\eta\co\left(W_1',T_A'\coprod T_B'\right)\to\left(W_2',T_A'\coprod T_B'\right)$$ satisfies the following condition on the boundary: $\eta(h_A)=h_A$, $\eta(h_B)=h_B$ and $\eta(\delta_A')=\delta_A'h_A^{u}$, $\eta(\delta_B')=\delta_B'h_B^{v}$. This implies that $W_1=W_1'/\l\delta_A'=1,\delta_B'=1\r$ is homeomorphic to  $W_2'/\l\delta_A'h_A^{u}=1,\delta_B'h_B^{v}=1\r$, where $W_i'/\l\delta_A'=1,\delta_B'=1\r$ denotes the space obtained from $W_i'$ after performing a Dehn filling along $T_A'\coprod T_B'$ identifying each curve $\delta_A'$ and $\delta_B'$ to the meridian of a solid torus.

On the other hand,   recall that it follows from our construction that $W_i=W_i'/\l\delta_A'=1,\delta_B'=1\r$.  Since $W_2'/\l\delta_A'=1,\delta_B'=1\r$ and   $W_2'/\l\delta_A'h_A^{u}=1,\delta_B'h_B^{v}=1\r$ are homeomorphic, this completes the proof of the lemma in Case 2.
 
  \emph{Case 3: $V_A$ is hyperbolic and $V_B$ is a Seifert fibered manifold.} Then by point (iii) of the lemma  we know that there exists a basis $\l\lambda_A,\mu_A\r$ of $\pi_1T_A$ such that $$(s_i)_{\ast}(\lambda_A)=h_B,\ \  (s_1)_{\ast}(\mu_A)=\delta_B\ \ {\rm and}\ \   (s_2)_{\ast}(\mu_A)=\delta_Bh^a_B$$ where $\delta_B$ denotes a section of $T_B$ with respect to the Seifert fibration of $V_B$, $h_B$ denotes the homotopy class of the generic fiber of $V_B$ and $a\in{\Z}$. Denote by $V_B'$, resp. $B'$, the space obtained from $V_B$, resp. $B$,  after removing a regular neighborhood $V(h_B)$ of a regular fiber in $V_B$. We write $\b V_B'=T_B\cup T'\cup T_1\cup...\cup T_p$ where $T'=\b V(h_B)$ and $T_1\cup...\cup T_p$ denote the components of $\b V_B\setminus T_B$.  For each component $T_i$, $i=1,...,p$ we fix a section $d_B^i$ of $T_i$ with respect to the Seifert fibration of $V_B'$ in such a way that $\pi_1V_B'$ admits a the following presentation:
 $$\l a_1,b_1,...,a_{g_B},b_{g_B},q_1,...,q_{r_B},\delta_B,\delta_B',d_B^1,...,d_B^p,h_B:$$ $$a_ih_Ba_i^{-1}=h_B, b_ih_Bb_i^{-1}=h_B, q_ih_Bq_i^{-1}=h_B, \delta_Bh_B\delta_B^{-1}=h_B, \delta'_Bh_B(\delta'_B)^{-1}=h_B,$$ 
 $$d_B^ih_B(d_B^i)^{-1}=h_B, q_j^{\gamma_j}=h_B^{\lambda_j}, [a_1,b_1]...[a_{g_B},b_{g_B}]q_1...q_{r_B}\delta_B\delta_B'.d_B^1...d_B^q=1\r$$ 
  where $\delta_B'$ is represented by the meridian of $V(h_B)$.
  
   We construct a homeomorphism $\eta\co V_B'\to V_B'$ such that 
$$\eta(\delta_B)=\delta_Bh_B^a, \eta(d^i_B)=d^i_B, \eta(a_i)=a_i, \eta(b_i)=b_i, \eta(q_i)=q_i\ \ {\rm and}\ \  \eta(\delta'_B)=\delta'_Bh_B^{-a}$$  
  
   Extend $\eta$ by the identity map on $A\coprod(B\setminus V_B)$. Denote by $W_i'$ the space obtained by gluing $A$ with $B'$ by identifying $T_A$ with $T_B$ via $s_i$. Then it follows from the construction that $W_1'$ is homeomorphic to $W_2'$.     Note that the homeomorphism $\eta\co(W_1',T')\to(W_2',T')$ satisfies the following condition on the boundary: $\eta(h_B)=h_B$ and $\eta(\delta')=\delta'_Bh_B^{-a}$. This implies that $W_1=W_1'/\l\delta_B'=1\r$ is homeomorphic to  $W_2'/\l\delta_B'h_B^{-a}=1\r$.
   
 On the other hand,     recall that it follows from our construction that $W_i=W_i'/\l\delta_B'=1\r$.  Since $W_2'/\l\delta_B'=1\r$ and   $W_2'/\l\delta_B'h_B^{-a}=1\r$ are homeomorphic, this completes the proof of the lemma in Case 3. The proof of Proposition \ref{C} is complete.
\end{proof}

  Combining  Proposition \ref{C} with Proposition \ref{B} we have completed the proof of Proposition \ref{A}. This ends the proof of Theorem \ref{dom}.

\section{Proof of Theorem \ref{poset}}

Let $M$ be a closed orientable graph manifold and let $(N_i)_{i\in{\N}}$ be a sequence of closed orientable Poincar\'e-Thurston 3-manifolds dominated by $M$ via degree one maps $f_i\co M\to N_i$. Since $M$ is a graph manifold then ${\rm Vol}(M)=0$ and since ${\rm Vol}(N_i)\leq{\rm Vol}(M)$ then the $N_i$'s are graph manifolds by Theorem \ref{inv}. Denote by $n_i$ the number of prime factors of each $N_i$ and write $N_i=N_i^1\sharp...\sharp N_i^{n_i}$ the prime decomposition of $N_i$ (see \cite{M}). First note that the sequence $(n_i)_{i\in{\N}}$ is finite. Indeed since $M$ 1-dominates $N_i$ then by point (i) of Proposition \ref{obvious} applied to degree one maps  we know that $r(M)\geq r(N_i)$, where $r(M)$ denotes the rank of the fundamental group of $M$. On the other hand, it follows from the Grushko's Theorem (\cite{Ma}) that $r(N_i)=r(N_i^1)+...+r(N_i^{n_i})$. This implies that the number of prime factors of the $N_i$'s is bounded by the rank of $\pi_1M$. Thus to complete the proof of Theorem \ref{poset} it is sufficient to show that the sequence of prime factors $(N_i^j)$, when $(i,j)\in{\N}\times\{1,...,n_i\}$ is finite up to homeomorphism.  For each $(i,j)\in{\N}\times\{1,...,n_i\}$, consider the projection map $p_i^j\co N_i\to N_i^j$. Such a map is a 0-pinch and has degree one. On the other hand the Cutting of Theorem of M. Gromov applied to the prime decomposition of each $N_i$ implies that ${\rm Vol}(N_i^j)={\rm Vol}(M)=0$ for all $i,j$. Then to complete the proof of Theorem \ref{poset} we apply Theorem \ref{dom} or   Theorem  \ref{known}   to the sequence of degree one maps defined by  $p_i^j\circ f_i$ depending on whether $N_i^j$ is Haken or not.

\end{document}